\documentclass[a4paper,11pt]{article}

\usepackage{amssymb,latexsym}
\usepackage{euscript}
\usepackage{eufrak}
\usepackage{amsmath}
\usepackage{mathrsfs}
\usepackage{tikz}
\usepackage{bussproofs}
\usepackage{bm}
\usepackage{enumerate}
\usepackage[ruled,vlined]{algorithm2e}
\usepackage{float}

\newtheorem{defi}{\bf Definition}[section]
\newtheorem{rem}{\bf Remark}[section]
\newtheorem{lem}{\bf Lemma}[section]
\newtheorem{prop}{\bf Proposition}[section]
\newtheorem{teo}{\bf Theorem}[section]
\newtheorem{cor}{\bf Corollary}[section]

\newtheorem{eje}{\bf Example}[section]

\def\n{\mathbf{n}}
\def\b{\mathbf{b}}
\def\ut{\mathbf{\frac 1 3}}
\def\dt{\mathbf{\frac 2 3}}
\def\0{\mathbf{0}}
\def\1{\mathbf{1}}

\def\Sup{\mathop{\rm Sup}\nolimits}
\def\Inf{\mathop{\rm Inf}\nolimits} 
\def\Hom{\mathop{\rm Hom}\nolimits}

\newcommand{\TSSix}{{\bf TS-Six}}
\newcommand{\GSix}{{\bf GSix}}
\newcommand{\SFSix}{{\bf SF-Six}}
\newcommand{\Six}{{\bf Six}}
\newcommand{\Sto}{{\bf S}}
\newcommand{\LFI}{{\bf LFI}}

\newcommand{\gsub}[1]{\mathop{{\sf gsub}}(#1)}

\newenvironment{dem}{\noindent\bf Proof. \rm}{\hfill $\mbox{\boldmath{$ \square$}}$}

\makeatother

\makeatletter

\newcommand*{\reflectmathsymbol}[2]{%
  \reflectbox{$\m@th#1#2$}%
}

\makeatletter
\providecommand*{\Dashv}{%
  \mathrel{%
    \mathpalette\@Dashv\models
  }%
}
\newcommand*{\@Dashv}[2]{%
  \reflectbox{$\m@th#1#2$}%
}
\makeatother

\title{ {\bf  Cut-free sequent-style systems for a logic associated to involutive Stone algebras}}

\author{
    Liliana M. Cant\'u 
    \and
    Mart\'{\i}n Figallo
}

\date{\textit{\small Departamento de Matem\'atica, Instituto de Matemática (INMABB), Universidad Nacional del Sur (UNS), Bah\'ia Blanca, Argentina}}

\begin{document}

\maketitle

\thispagestyle{myheadings}

\begin{abstract}

In \cite{LC, LCMF}, it was introduced a logic (called \Six ) associated to a class of algebraic structures known as {\em involutive Stone algebras}. This class of algebras, denoted by \Sto ,  was considered by the first time in \cite{CS1} as a tool  for the study of certain problems connected to the theory of finite--valued \L ukasiewicz--Moisil algebras. In fact, \Six\ is the logic that preserves degrees of truth with respect to the class \Sto .

Among other things, it was proved that \Six\ is a 6--valued logic that is a {\em Logic of Formal Inconsistency} (\LFI ); moreover, it is possible to define a consistency operator in terms of the original set of connectives. A Gentzen-style system (which does not enjoy the cut--elimination property) for \Six\ was given. Besides, in \cite{LCMF}, it was shown that \Six\ is matrix logic that is determined by a finite number of matrices; more precisely, four matrices. However, this result was further sharpened in \cite{SMUR} proving that just one of these four matrices determine \Six , that is, \Six\ can be determined by a single 6--element logic matrix.  

In this work, taking advantage of this last result, we apply a method due to Avron, Ben-Naim and Konikowska \cite{Avron02} to present different Gentzen systems for \Six\ enjoying the cut--elimination property. This allows us to draw some  conclusions about the method as well as to propose additional tools for the streamlining process. Finally, we present a decision procedure for \Six\ based in one of these systems.

\vspace*{4mm}


\vspace*{4mm}

\end{abstract}

\section{\large \bf Introduction}\label{s1}

Pursuing to characterize five-valued \L ukasiewicz-Moisil algebras as De Morgan algebras such that their complemented elements satisfy certain conditions; and to obtain the corresponding characterizations for the four and three-valued cases, R. Cignoli and M. Sagastume introduced the class of {\em involutive Stone algebras} and denoted it by \Sto , in \cite{CS1}.
These same authors proved (among other things) that \Sto\ is a variety, provided a topological duality for them and used this duality for determining the subdirectly irreducible and simple objects of this class. From an algebraic point of view, \Sto -algebras are De Morgan algebras endowed with an additional unary operation,  denoted by $\nabla$. More precisely, let $A$ be a De Morgan algebra and let $K(A)$ be the set of all elements $a\in A$ such that the De Morgan negation of $a$, $\neg a$, coincides with the Boolean complement of $a$. For every $a\in A$, let $K_a=\{k\in K(A) : a\leq k\}$ and, if $K_a$ has a least element, denote it by $\nabla a$. Then, the class of involutive Stone algebras is the class of all De Morgan algebras $A$ such that for every $a\in A$, $\nabla a$ exists, and the map $a \mapsto \nabla a$ is a lattice--homomorphism.

In \cite{LC, LCMF}, it was considered a logic associated to the class of involutive  Stone algebras by the first time. More precisely, it was defined the logic that preserves degrees of truth associated to involutive Stone algebras, named \Six . This follows a very general pattern that can be considered for any class of truth structure endowed with an ordering relation; and which intend to exploit many-valuedness focusing on the notion of inference that results from preserving lower bounds of truth values, and hence not only preserving the greatest value $1$. 
Among other things, it was proved that \Six\  is a many valued logic (with six truth values) that can be determined by a finite number of matrices (four matrices). Besides, \Six\ turns out to be  a paraconsistent logic, moreover, it was proved that it is a genuine {\bf LFI} (Logic of Formal Inconsistence) with a consistency operator that can be defined in terms of the original set of connectives. Furthermore, a Gentzen-style system for \Six\ was given in \cite{LC, LCMF}. Unfortunately, this system does not enjoy the cut--elimination property (as we shall see below).

Six-valued logics have emerged in the past for diverse reasons and in different contexts. In \cite{GaMo}, the authors introduce a six-valued logic useful in representing incomplete knowledge. A practical advantage of this logic is that it allows a system to progressively reason about what it will or will not know (or what can or cannot happen) as time advances and further knowledge is acquired from the external world.  To capture
the changing nature of the world, they make a distinction between two types of unknown knowledge values:
permanently or eternally unknown values, and temporarily unknown values. Remarkably, one of the two orders proposed for this six-element set of truth values coincides with the order of $\mathbb{S}_6$ of Figure \ref{fig1}.
Other examples of six-valued logics can be found in \cite{FoMo, GaLoPaFeDi, EnTu}.

Recently, the logics that arise from the involutive Stone algebras has gained some attention. In \cite{SMUR}, it was proved that \Six\ is a conservative expansion of the Belnap-Dunn four-valued logic (i.e. the order-preserving logic of the variety of De Morgan algebras), and gave a finite Hilbert-style axiomatization for it, from which a decision procedure is available. More generally, the authors introduced a method for expanding conservatively every super-Belnap logic so as to obtain an extension of \Six ; and showed that every logic thus defined can be
axiomatized by adding a fixed finite set of rule schemata to the corresponding super Belnap base logic.
On the other hand, in \cite{JGVGSMJMUR}, it was introduced the class of {\em perfect paradefinite algebras} which are De Morgan algebras enriched by a {\em perfection operator} (see \cite{JGVGSMJMUR}). These algebras turn out to be term-equivalent to involutive Stone algebras. Later, these authors studied the $1$-assertional logic and also the order-preserving logic associated to perfect paradefinite algebras.
For all this, we think that \Six\ is an interesting example of $6$-valued logic with a nice algebraic semantics that deserves further study.  

\

On the other hand, in \cite{Avron01}, the authors use the Rasiowa-Sikorski  decomposition
methodology to get sound and complete proof systems employing the formalism of $n$-sequents
for all propositional logics based on non-deterministic matrices. Later, in \cite{Avron02}, a general method to transform a given sound and complete $n$-sequents proof systems into an equivalent  sound and complete system of ordinary two-sided sequents, for languages satisfying a certain minimal expressiveness condition, was presented. The systems obtained in these steps are hardly optimal respecting both the number and size of rules (see \cite{Avron01, Avron02}). As we shall see, in this work we obtain two proof systems for \Six\, one in terms of $6$-sequents and the other using ordinary two-sided sequents, containing 84 rules the former and 230 rules the latter. Besides, many of these rules have up to six upper sequents to introduce a binary connective. That is the reason why three general streamlining principles, to reduce them to a more compact form, were also proposed in   \cite{Avron01, Avron02}. However, the streamlinig process is not algorithmical at all and strongly depends on the choices made in the different stages and the skills and abilities of whom is leading this process.

\

The focus of this paper is two-fold; the first objective is to present a sequent system for \Six\ enjoying the cut elimination property. As it is well-known, the cut-elimination theorem plays a mayor role in proof theory since it allows to have control on the size and structure of proofs; and therefore, many properties of the logic can be derived analyzing the structure of proofs (as for example the various interpolation theorems).  For this, we take advantage of a result in \cite{SMUR}, which state that \Six\ can be determined by a single 6--element matrix.  Then, following the abovementioned method of Avron et al., we provide new Gentzen-style systems for \Six\; all of them enjoying the cut--elimination property. The second objective is, using the experience gained by achieving our first objective, to propose new principles and tools for the streamlining process pursuing some kind of ``mechanization'' for some parts of it. This allows us to obtain a cut-free sequent system with 26 rules, which introduce (interaction of) connectives, with one or two upper sequents. Finally, we explicit a decision procedure for \Six\ based in this last system which is an adaptation of the original decision process given by Gentzen (see \cite{Tau}) for classical and intuitionistic propositional logic; and from which we obtain bottom-up proof search.

\

\section{\large \bf Preliminaries}\label{s2}

Recall that a {\em De Morgan algebra} is a structure $\langle A, \wedge, \vee, \neg, 0, 1 \rangle$ of type $(2,2,1,0,0)$ such that the reduct $\langle A, \wedge, \vee,  0, 1 \rangle$ is a bounded distributive lattice and $\neg$ satisfies the identities:  
$$\neg \neg x \approx x, \mbox{ \, and }$$  
$$\neg(x \wedge y) \approx \neg x \vee \neg y.$$
Let $A$ be a De Morgan algebra, $a \in A$ is said to be {\em complemented} if there exists $b\in A$, such that $a \vee b = 1$ and $a \wedge b = 0$; it is said that $b$ is the {\em Boolean complement} of $a$ and it is denoted $a'$. If the complement of an element $a$ exists, it is unique. We denote by $B(A)$ the {\em center of} $A$, that is, the set of all complemented elements of $A$.  It can be seen that, $B(A)$ is a  De Morgan subalgebra of $A$. Let $K(A)$ be the subset of $A$ formed by all elements $a \in B(A)$ such that the De Morgan negation of $a$, $\neg a$, coincides with its complement $a'$. For every $a\in A$, let $K_a =\{k \in K(A): a\leq k \}$. If $K_a$ has a least element, we denote it by $\nabla a$ as in \cite{CS1}. If for every $a\in A$, $K_a$ has a least element, it is said that $A$ is a $\nabla$-De Morgan algebra. The {\em class of involutive Stone algebras}, denoted by {\bf S}, is the subclass of $\nabla$-De Morgan algebras with the property that the map $a \mapsto \nabla a$ from $A$ into $K(A)$ is a lattice-homomorphism. \\
In \cite{CS1}, it was proved that $\bf S$ is an equational class. Indeed, a De Morgan algebra  $A \in {\bf S}$ if and only if there is an operator $\nabla : A\to A$  satisfying the following equations:
$$\nabla 0 \approx 0,$$
$$a \wedge \nabla a \approx a,$$
$$\nabla (a \wedge b) \approx \nabla a \wedge \nabla b,$$
$$\neg \nabla a \wedge \nabla a \approx 0.$$

\noindent For each $n\geq 2$, {\bf \L}$_{n}$ will denote the $n$-element \L ukasiewicz  chain \, $0 < \frac{1}{n-1} < \dots < \frac{n-2}{n-1} < 1$ with the well-known lattice structure and where \, $\neg (\frac{i}{n-1}) = 1 - \frac{i}{n-1}$ for $i=0,1,\dots,n-1$. The algebras {\bf \L}$_{n}$ are  examples of involutive Stone algebras, defining $\nabla x= 1$, for all $x\not= 0$ and $\nabla 0=0$.
The subdirectly irreducible algebras are {\bf \L}$_n$, for $2\leq n\leq5$ and the {\bf S}-algebra  $\mathbb{S}_6$ given in Figure \ref{fig1}.

\begin{figure}[h]
\centering
\begin{center}
\begin{tikzpicture}[scale=.7]
  \node (one) at (0,2) {$\1$};

  \node (zero) at (0,0) {$\0$};
 \node () at (0,-2) {{\bf \L}$_2$};
  \draw (zero) -- (one);

\end{tikzpicture}  \hspace{2cm}
\begin{tikzpicture}[scale=.7]
  \node (one) at (0,4) {$\1$};
  \node (a) at (0,2) {$\bf \frac{1}{2}$};
  \node (zero) at (0,0) {$\0$};
 \node () at (0,-2) {{\bf \L}$_3$};
  \draw (zero) -- (a) -- (one);

\end{tikzpicture}  \hspace{2cm}
\begin{tikzpicture}[scale=.7]
  \node (one) at (0,6) {$\1$};
  \node (b) at (0,4) {$\dt$};
  \node (a) at (0,2) {$\ut$};
  \node (zero) at (0,0) {$\0$};
\node () at (0,-2) {{\bf \L}$_4$};
  \draw (zero) -- (a) -- (b) -- (one);
\end{tikzpicture} \hspace{2cm}
\begin{tikzpicture}[scale=.7]
  \node (one) at (0,8) {$\1$};
  \node (b) at (0,6) {$\bf \frac{3}{4}$};
  \node (c) at (0,4) {$\bf \frac{2}{4}$};
  \node (a) at (0,2) {$\bf \frac{1}{4}$};
  \node (zero) at (0,0) {$\0$};
\node () at (0,-2) {{\bf \L}$_5$};
  \draw (zero) -- (a) -- (c) -- (b) -- (one);
\end{tikzpicture} \hspace{2cm}
\begin{tikzpicture}[scale=.7]
 	
	\node (one) at (0,6)  {$\1$};
	\node (c) at (0,4.5) {$\dt$};
  \node (b) at (-1.5,3)  {$\n$};
	 \node (n) at (1.5,3)  {$\b$};
 \node (a) at (0,1.5) {$\ut$};
\node (zero) at (0,0) {$\0$};
\node () at (0,-2) {$\mathbb{S}_6$};
  \draw (zero) -- (a) -- (n) -- (c) -- (one) -- (c) -- (b) -- (a);

\end{tikzpicture} 
\end{center}
\noindent where the operations in $\mathbb{S}_6$ are given by the diagram above and \, $\neg \0 =\1$, $\neg \1 =\0$, $\neg \n =\n$, $\neg \b =\b$, $\neg \ut =\dt$ and $\neg \dt=\ut$. Besides,  $\nabla x= \1$, for all $x\not= \0$ and $\nabla \0=\0$. \\
\caption{Subdirectly irreducible algebras of $\bf S$.}
\label{fig1}
\end{figure}

\noindent From the fact that {\bf \L}$_n$, for $2\leq n\leq5$, are $\bf S$-subalgebras of $\mathbb{S}_6$, we have the following:

\begin{prop}\label{coro1}(\cite{CS1}) $\mathbb{S}_6$ generates the variety {\bf S}.
\end{prop}

In \cite{LC, LCMF}, it was studied the logic that preserves degrees of truth associated to involutive Stone algebras. This follows a very general pattern that can be considered for any class of truth structure endowed with an ordering relation; and which intends to exploit manyvaluedness focusing on the notion of inference that results from preserving lower bounds of truth values, and hence not only preserving the value $1$. More precisely, let $\mathfrak{Fm}=\langle Fm, \wedge, \vee, \neg, \nabla, \bot, \top \rangle$ be the absolutely free algebra of type $(2,2,1,1,0,0)$ generated by some fixed denumerable set of propositional variables $Var$. As usual, the letters $p,q, \dots$ denote propositional variables, the letters $\alpha,\beta, \dots$ denote formulas and $\Gamma,\Sigma, \dots$ 
sets of formulas. Let $\Gamma\cup\{\alpha\}\subseteq Fm$, then \, $\Gamma \models_{\mbox{\Six}} \alpha$ \, iff \, there exists $\Gamma_0\subseteq \Gamma$ finite such that  $\forall A\in {\bf S}$, $\forall v\in \Hom(\mathfrak{Fm},A)$, $\forall a\in A$,  if $v(\gamma)\geq a$, for all $\gamma\in \Gamma_0$, then $v(\alpha)\geq a$. Besides,  $\emptyset \models_{\mbox{\Six}} \alpha$ \,  iff \, $\forall A\in {\bf S}$, $\forall v\in \Hom(\mathfrak{Fm},A)$, $v(\alpha)=1$. \\
It is clear that \Six\ is a sentential logic, that is, $\models_{\mbox{\Six}}$ is a finitary consequence relation over $Fm$. Besides,

\begin{prop}(\cite{LC, LCMF}) Let $\{\alpha_1,\dots, \alpha_n, \alpha\}\subseteq Fm$, $n\geq 1$. Then, \,  $\alpha_1,\dots, \alpha_n \models_{\mbox{\Six}} \alpha$ \, iff \,  the inequality \,  $\alpha_1 \wedge\dots \wedge \alpha_n  \preccurlyeq \alpha$ holds in the variety {\bf S}.
\end{prop}

\noindent From this,  we have

\begin{prop}(\cite{LC, LCMF})  Let $\alpha,\beta \in Fm$. Then,
\begin{itemize}
\item[(i)] $\alpha \models^{\leq}_{\bf S} \beta$ \, iff \, the inequality \, $\alpha \preccurlyeq \beta $, holds in {\bf S},
\item[(ii)] $\alpha \Dashv \models^{\leq}_{\bf S} \beta$ \,  iff  \, the equation \, $\alpha \approx \beta $ holds in {\bf S}.
\end{itemize}
\end{prop}

\

\noindent It was shown in \cite{LCMF}, that \Six\  is a matrix logic determined by four matrices.

\begin{prop} The logic that preserves degrees of truth associated to involutive Stone algebras, \Six , is a six--valued logic determined by the matrices $\langle \mathbb{S}_6, \left[\ut\right)\rangle$, $\langle \mathbb{S}_6, [\b)\rangle$, $\langle \mathbb{S}_6, \left[\dt\right)\rangle$ and $\langle \mathbb{S}_6, [\1)\rangle$.
\end{prop}

\noindent Here, $[a)$ represents the lattice filter of $\mathbb{S}_6$ generated by $a$. Later this result was further sharpened in \cite{SMUR}, showing that \Six\ can be determined by a single logic matrix.

\begin{prop} (\cite{SMUR}) \,  \Six\ is the logic determined by the matrix $\langle \mathbb{S}_6, [\b)\rangle$.
\end{prop}

\begin{rem} The well--known Dunn-Belnap's four--element given by the diagram

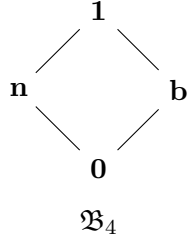
\begin{center}
\begin{tikzpicture}[scale=.7]

	\node (one) at (-1,3)  {$\1$};
  \node (b) at (-2.5,1.5)  {$\n$};
	\node (a) at (0.5,1.5) {$\b$};
  \node (zero) at (-1,0)  {$\0$};
\node () at (-1,-1) {$\mathfrak{B}_4$};
  \draw (zero) -- (a) -- (one) -- (b) -- (zero);
	\end{tikzpicture}
\end{center}

\noindent with the negation defined by $\neg \0=\1$, $\neg \n=\n$, $\neg \b=\b$ and $\neg \1=\1$, on the one hand, and {\bf \L}$_4$, in the other, are De Morgan subalgebras of $\mathbb{S}_6$, respectively. This fact explains our motivations for the way in which we named the elements of $\mathbb{S}_6$.
\end{rem}

In \cite{LCMF}, it was given a multiple-conclusioned  sequent system for \Six . Let $\mathfrak{S}$  be the  sequent calculus whose  axioms  and  rules  are  the following:

\

\noindent {\bf Axioms}
$$ \mbox{(structural axiom) \, } \displaystyle {\alpha \Rightarrow \alpha} \hspace{2cm} \mbox{($\bot$) \,  } {\bot\Rightarrow } \hspace{2cm} \mbox{($\top$) \,  } {\Rightarrow \top }$$
$$\mbox{(First modal axiom) \,  } {\alpha\Rightarrow \nabla\alpha} \hspace{2cm} \mbox{(Second modal axiom) \,  } {\Rightarrow \nabla\alpha \vee \neg \nabla \alpha}$$

\noindent {\bf Structural rules}
$$ \mbox{(Left weakening) \, } \displaystyle \frac{\Gamma \Rightarrow \Sigma} {\Gamma, \alpha \Rightarrow \Sigma} \hspace{2cm} \mbox{(Right weakening) \, }
 \displaystyle \frac{\Gamma \Rightarrow \Sigma} {\Gamma \Rightarrow \Sigma, \alpha} $$
$$ \mbox{(Cut) \, }
 \displaystyle \frac{\Gamma \Rightarrow \Sigma, \alpha  \hspace{0.5cm} \alpha, \Gamma \Rightarrow \Sigma}{\Gamma \Rightarrow \Sigma} $$

\noindent {\bf Logic Rules}
$$ \mbox{($\wedge \Rightarrow$) \, } \displaystyle \frac{\Gamma, \alpha, \beta \Rightarrow \Sigma} {\Gamma, \alpha \wedge \beta \Rightarrow \Sigma} \hspace{2cm} \mbox{($\Rightarrow \wedge$) \, }
 \displaystyle \frac{\Gamma \Rightarrow \Sigma, \alpha  \hspace{0.5cm} \Gamma \Rightarrow \Sigma, \beta}{\Gamma \Rightarrow \Sigma, \alpha \wedge \beta} $$
$$ \mbox{($\vee \Rightarrow$) \, } \displaystyle \frac{\Gamma, \alpha \Rightarrow \Sigma \hspace{0.5cm} \Gamma, \beta \Rightarrow \Sigma} {\Gamma, \alpha \vee \beta \Rightarrow \Sigma} \hspace{2cm}  \mbox{($\Rightarrow \vee$) \, }
 \displaystyle \frac{\Gamma \Rightarrow \Sigma, \alpha, \beta }{\Gamma \Rightarrow \Sigma, \alpha \vee \beta} $$
$$ \mbox{($\neg$) \, } \displaystyle \frac{\alpha \Rightarrow \beta } {\neg \beta \Rightarrow \neg \alpha } \hspace{2cm} \mbox{($\neg \neg \Rightarrow$) \, } \displaystyle \frac{\Gamma, \alpha \Rightarrow \Sigma } {\Gamma, \neg \neg \alpha \Rightarrow \Sigma} \hspace{2cm} \mbox{($\Rightarrow \neg \neg$)} \, \frac{\Gamma \Rightarrow \alpha, \Sigma } {\Gamma \Rightarrow \neg \neg \alpha, \Sigma} $$
$$ \mbox{($\nabla$) \, } \displaystyle \frac{\Gamma, \alpha \Rightarrow \nabla \Sigma } {\Gamma, \nabla \alpha \Rightarrow \nabla \Sigma} \hspace{2cm} \mbox{($\neg \nabla \Rightarrow$) \, } \displaystyle \frac{\Gamma, \neg \nabla \alpha \Rightarrow \Sigma } {\Gamma, \nabla \neg \nabla \alpha \Rightarrow \Sigma}   $$

\

\

Unfortunately, the calculus $\mathfrak{S}$ does not enjoy the cut-elimination property. Indeed, it is not difficult to check that, if $p$ and $q$ are two propositional variables then the sequent 
$$ p\vee q \Rightarrow \neg(\neg p\wedge \neg q)$$ 
\noindent is provable in $\mathfrak{S}$ and every proof of it use necessarily the cut rule. 

\

On the other hand, \Six\  is a paraconsistent logic in the sense of da Costa (cf.~\cite{daC,daC2}) since inferences in the logic not necessarily are trivialized in the presence of contradictions.

\begin{prop} (\cite{LC, LCMF}) {\bf \em Six} is non--trivial, non--explosive and  paracomplete.
\end{prop}

\noindent Besides, in the language of  \Six\  let $\bigcirc(p)=\{ \neg \nabla(p \wedge  \neg p) \}$. Then,

\begin{lem} (\cite{LC, LCMF}) {\bf \em Six} is finitely gently explosive with respect to $\bigcirc(p)$ and $\neg$.
\end{lem}

\noindent From all the above we have:

\begin{teo} (\cite{LC, LCMF}) {\bf \em Six} is an {\bf LFI} with respect to $\neg$ and with consistency operator $\circ$ defined by $\circ \alpha =  \neg \nabla(\alpha \wedge  \neg \alpha) $, for all $\alpha \in Fm$.
\end{teo}

\noindent Moreover, \ $\models_{\mbox{\small \bf \em Six}} \circ \neg^{n} {\circ}
\alpha$, for all $n\geq 0$. In particular, $\models_{\mbox{\small \bf \em Six}}
{\circ} {\circ} \alpha$. So, {\bf \em Six} validates all axioms $(cc)_{n}$ of the logic {\bf mCi} (see \cite{WCMCJM}).
As it is usual in the framework of {\bf LFI}'s, it is also possible to define an inconsistency operator~$\bullet$ \, in {\bf \em Six} \, in the following way:

$$ \bullet \alpha =_{def} \neg \circ \alpha.$$

\

\noindent Then, $\bullet \alpha$ is logically equivalent to $\nabla(\alpha \wedge \neg \alpha)$. Note that $\bullet$ and $\circ$ can be expressed equivalently as  \, $\bullet \alpha = \nabla( \alpha \wedge \neg \alpha)$ \, and \, $\circ \alpha
= \Delta(\alpha \vee \neg \alpha)$. The truth--table of both connectives is the following:

\

\begin{center}
\begin{tabular}{|c|c|c|}         \hline
$p$ & $\circ p$ & $\bullet p$\\ \hline
$\0$       & $\1$  & $\0$         \\ \hline
$\ut$ & $\0$ & $\1$  \\ \hline
$\n$       & $\0$ & $\1$           \\ \hline
$\b$       & $\0$ & $\1$          \\ \hline
$\dt$ & $\0$ & $\1$  \\ \hline
$\1$       & $\1$  &$\0$          \\ \hline
\end{tabular}
\end{center}
\

Finally, in what follows, we shall present the matrix $\langle \mathbb{S}_6, [\b)\rangle$ in the context of \cite{Avron02}. More precisely, let ${\cal T}_6=\{\0,\ut, \n, \b, \dt, \1\}$ be the set of truth-values with an underlying order given by the diagram in Figure \ref{fig1}; then let   ${\cal M}_6$ be the matrix $\langle {\cal T}_6, {\cal D}, {\cal O} \rangle$  where ${\cal D}=\{\b, \dt, \1\}$  is the set of designated values and ${\cal O}=\{\hat{\wedge},\hat{\vee},\hat{\neg},\hat{\nabla}\}$ where, for $x, y \in {\cal T}_6$,
$x\hat{\wedge} y= \Inf \{x,y\}$ \,  $x \hat{\vee} y= \Sup \{x,y\}$, and the operations $\hat{\neg}$ and $\hat{\nabla}$ are given by the tables

\begin{center}
\begin{tabular}{|c|c|c|} \hline 
$x$ & $\hat{\neg} x$ & $\hat{\nabla} x$ \\ \hline \hline
$\0$ & $\1$ & $\0$ \\ 
$\ut$ & $\dt$ & $\1$ \\ 
$\n$ & $\n$ & $\1$ \\ 
$\b$ & $\b$ & $\1$ \\ 
$\dt$ & $\ut$ & $\1$ \\ 
$\1$ & $\0$ & $\1$ \\ \hline
\end{tabular}
\end{center}

\

\section{Cut-free sequential systems for matrix logics} \label{s3}

In \cite{Avron01}, the authors use the Rasiowa-Sikorski  decomposition
methodology to get sound and complete proof systems employing $n$-sequents
for all propositional logics based on non-deterministic matrices. Later, in \cite{Avron02}, it was presented a general method to transform a given sound and complete $n$-sequents proof systems into an equivalent  sound and complete system of ordinary two-sided sequents for languages satisfying a certain minimal expressiveness condition. In this section  we shall recall both methods considering ordinary (deterministic) matrices. 

\

Let $\mathscr{L}$ be a propositional language and let $\mathfrak{Fm}$ be the absolutely free algebra over $\mathscr{L}$ generated by some denumerable set of variables, with underlying set (of formulas) $Fm$.
Recall that an ($n$-valued) matrix  ${\cal M}$  for $\mathscr{L}$  is a triple $\langle {\cal T}, {\cal D }, {\cal O}\rangle$  where  ${\cal T}$ is a finite, non-empty set of ($n$) truth values, ${\cal D}$ is a non-empty proper set of ${\cal T}$, and ${\cal O}$ includes a $k$-ary function $\hat{f}: {\cal T}^k\to {\cal T}$ for each $k$-ary connective $f$ in the propositional language.\\[2mm]
Recall that a valuation in ${\cal M}$ is a function $v:Fm\to {\cal T}$ such that
$$v(f(\psi_1, \dots, \psi_k))=\hat{f}(v(\psi_1), \dots, v(\psi_k))$$
for each $k$-ary connective $f$ and all $\psi_1, \dots, \psi_k\in Fm$. 
A formula $\alpha \in Fm$ is satisfied by a given valuation $v$, in symbols $v\models \alpha$, if $v(\alpha)\in {\cal D}$. A sequent $\Gamma \Rightarrow \Delta$ is satisfied by the valuation $v$, in symbols $v\models \, \Gamma \Rightarrow \Delta$, if either $v$ does not satisfy some formula in $\Gamma$ or $v$ satisfies some formula in $\Delta$. A sequent is {\em valid} (w.r.t the matrix ${\cal M}$) if it is satisfied by all valuations. We write $\vdash_{\cal M}  \Gamma\Rightarrow \Delta$ to indicate that the sequent  $\Gamma\Rightarrow \Delta$ is valid in $\cal M$.\\ 
Now, suppose that ${\cal T}=\{t_0,\dots, t_{n-1}\}$, where $n\geq 2$,  and ${\cal D}=\{t_d,\dots, t_{n-1}\}$, where $1\leq d \leq n-1$.\\
Intuitively, $n$-sequents consist of $n$ sets of formulas, one for each truth value and, as we shall see, a valuation $v$ validates a given $n$-sequent if there is a formula in the $i$-th set to which $v$ assigns the $i$-th truth value, for some $i$, $0\leq i\leq n-1$ .

\begin{defi} (see \cite{Avron01}) An $n$--sequent over $\mathscr{L}$ is an expression  
$$ \Gamma_0 \mid \dots \mid \Gamma_{n-1}$$
where, for each $i$, $\Gamma_i$ is a finite set of formulas. A valuation $v$ satisfies the $n$--sequent $\Gamma_0 \mid \dots \mid \Gamma_{n-1}$ iff there exists $i$, $0\leq i \leq n-1$ \, and $\psi\in \Gamma_i$ such that $v(\psi)=t_i$. An $n$--sequent is valid if it is satisfied by every valuation $v$. 
\end{defi}  

\noindent Note that, a valuation $v$ satisfies an ordinary sequent $\Gamma  \Rightarrow \Delta$ iff $v$ satisfies the $n$--sequent $\Gamma_1 \mid \dots \mid \Gamma_{n-1}$ where $\Gamma_i=\Gamma$ for all $0\leq i\leq d-1$ and $\Gamma_j=\Delta$ for all $d\leq j\leq n-1$ .

\noindent An alternative presentation of $n$-sequents is by means of sets of signed formulas. A signed formula over the language $\mathscr{L}$ and ${\cal T}$, is an expression of the form
$$t_i : \psi$$
where $t_i\in  {\cal T}$ and $\psi \in Fm$. A valuation $v$ satisfies the signed formula $t_i : \psi$ iff $v(\psi)=t_i$. If $\Omega\subseteq {\cal T}$ and $\Gamma \subseteq Fm$, we denote by $\Omega : \Gamma$ the set
$$\Omega : \Gamma =\{ t:\alpha \mid t\in \Omega, \alpha \in \Gamma\}$$
If $\Omega=\{t\}$, we write $t : \Gamma$ instead of $\{t\} : \Gamma$.  
A valuation satisfies the set of signed formulas $\Omega : \Gamma$  if it satisfies some signed formula of $\Omega : \Gamma$; and we say that $\Omega : \Gamma$ is valid if it is satisfied by every valuation $v\in {\cal V}$.
It is clear that,  the $n$--sequent $ \Gamma_0 \mid \dots \mid \Gamma_{n-1}$ is valid iff the set of signed formulas $\bigcup \limits_{i=0}^ {n-1} t_i : \Gamma_i$ is valid.

\

\noindent In \cite{Avron01} it was developed a generic $n$-sequent system for any logic based on an $n$-valued matrix. Consider the $n$-valued matrix ${\cal M}=\langle {\cal T}, {\cal D }, {\cal O}\rangle$ and let $SF_{\cal M}$ the system defined as follows:

\begin{itemize}
\item {\bf Axioms:} \, \,    ${\cal T}: \alpha$
\item {\bf Structural rules: } weakening
\item {\bf Logical rules: } for each $k$-ary connective $f$ and every $(a_1, \dots, a_k) \in   {\cal T}^k$ 

$$ \displaystyle  \frac{\Omega, a_1: \alpha_1 \, \, \dots \, \, \Omega, a_k: \alpha_k } {\Omega,  \hat{f}(a_1, \dots, a_k) : f(\alpha_1, \dots,\alpha_k)}$$
\end{itemize}

\

\begin{teo}\label{TeoAv01} (\cite{Avron01}) The system $SF_{\cal M}$ is sound and complete w.r.t. the matrix ${\cal M}$
\end{teo}

It is possible to translate a given $n$-sequent calculus of an $n$-valued logic satisfying certain expressiveness conditions to an ordinary two-sided sequent calculus (see \cite{Avron02}). Let $Fm_p$ be the set of all formulas of  $Fm$ that have $p$ as their only propositional variable, i.e., $Fm_p=\{ \alpha  \in Fm : Var(\alpha)=\{p\} \}$. Let us denote by ${\cal N}$ the set ${\cal T}\setminus {\cal D}$.

\begin{defi}(\cite{Avron02}) \label{SufficientlyExp} The language $\mathscr{L}$ is sufficiently expressive for ${\cal M}$ iff for any $i$, $0\leq i \leq n-1$ there exist natural numbers $l_i,m_i$  and formulas $\alpha_{j}^{i}, \beta_{k}^{i}\in Fm_p$, for $1\leq j\leq l_i$ and $1\leq k\leq m_i$ such that for any valuation $v$, the following conditions hold:\\
(i) $\alpha_{1}^{i}=p$ if $t_i\in {\cal N}$ \, and \, $\beta_{1}^{i}=p$ if $t_i\in {\cal D}$, \\
(ii) For any $\varphi \in Fm$ and $t_i\in {\cal T}$
$$v(\varphi)=t_i \, \Leftrightarrow \ v(\alpha_{1}^{i}[p/\varphi]), \dots, v(\alpha_{l_i}^{i}[p/\varphi]) \in {\cal N} \, and \,  v(\beta_{1}^{i}[p/\varphi]), \dots, v(\beta_{m_i}^{i}[p/\varphi]) \in {\cal D}$$
where $\alpha_{j}^{i}[p/\varphi]$ ($\beta_{k}^{i}[p/\varphi]$) is the formula obtained by the substitution of $p$ by $\varphi$ in $\alpha_{j}^{i}$ ($\beta_{k}^{i}$).
\end{defi}  

\noindent If $\Gamma$ is a set of formulas and $\alpha\in Fm_p$, we denote by $\alpha[p/\Gamma]$ the set

$$\alpha[\Gamma]=\{\alpha[p/\gamma] \mid \gamma\in \Gamma\}$$
The translation process consists of replacing each $n$-sequent by a semantically equivalent set of two-sided sequents. \\
Suppose that $\mathscr{L}$ is a sufficiently expressive language and for $0\leq i \leq n-1$ let $l_i$, $m_i$, $\alpha_{j}^{i}$  and $\beta_{k}^{i}$ as in Definition \ref{SufficientlyExp}. Consider the $n$--sequent $\Sigma = \Gamma_0 \mid \dots \mid \Gamma_{n-1}$ over $\mathscr{L}$. A partition $\pi$ of the $n$--sequent $\Sigma$ is a tuple $\pi=(\pi_0,\dots,\pi_{n-1})$  such that, for every $i$, $\pi_i$ is a partition of the set $ \Gamma_i$ of the form:
$$\pi_i=\{\Gamma'_{ij}\mid 1\leq j\leq l_i\}\cup\{\Gamma''_{ik}\mid 1\leq k\leq m_i\}$$

\

\begin{eje} Suppose that ${\cal T}=\{t_0, t_1, t_2\}$, that is,  $n=3$   and  $0\leq i\leq 2$. Let $l_0=2=l_2$, $l_1=1$, $m_0=1$, $m_1=2$ and $m_2=0$ (see Definition \ref{SufficientlyExp}). Then,  a partition $\pi$ of a given $3$-sequent is of the form

$$\big( \big\{ \Gamma'_{01}, \Gamma'_{02}, \Gamma''_{01}\big\},  \big\{ \Gamma'_{11}, \Gamma''_{11}, \Gamma''_{12}\big\}, \big\{ \Gamma'_{21}, \Gamma'_{21}\big\}\big)$$ 
Note that, in this context, empty cells are allowed. So, if we have the $3$-sequent \, $\phi, \varphi \, | \, \delta \,  | \, \gamma, \rho$ \break the next is a partition of it.    
$$\big( \big\{ \{\phi\}, \emptyset, \{\varphi\}\big\},  \big\{ \emptyset, \{\delta\}, \emptyset \big\}, \big\{ \{\gamma, \rho\}, \emptyset\big\}\big)$$ 
\end{eje}

\

\noindent Then, given a partition $\pi$ of the $n$-sequent $\Sigma$, we define the two-sided sequent $\Sigma_{\pi}$ determined by $\Sigma$ and the partition $\pi$,  as follows:
$$\bigcup \limits_{j=1}^{l_0} \alpha_{j}^{0}[\Gamma'_{0j}], \dots, \bigcup \limits_{j=1}^{l_{n-1}} \alpha_{j}^{n-1}[\Gamma'_{(n-1)j}] \, \Rightarrow \, \bigcup \limits_{k=1}^{m_0} \beta_{k}^{0}[\Gamma''_{0k}], \dots, \bigcup \limits_{k=1}^{m_{n-1}} \beta_{k}^{n-1}[\Gamma''_{(n-1)k}]  $$
Let $\Pi$ be the set of all partitions of the $n$--sequent $\Sigma$, it is denoted by $TWO(\Sigma)$ the set 
$$TWO(\Sigma)= \{\Sigma_{\pi}\mid \pi\in \Pi\}$$
Then, if ${\cal C}$ is an  $n$--sequent calculus over $\mathscr{L}$. We denote   $TWO({\cal C})$ the (ordinary) sequent calculus over $\mathscr{L}$ given by:
\begin{itemize}
\item[] {\bf Axioms:} $TWO(A)$, for all axiom $A$ of $\cal C$,
\item[] {\bf Inference rules:} $\displaystyle \frac{TWO(S)}{\Sigma'}$, where $\displaystyle \frac{S}{R}$ is a rule in $\cal C$ and $\Sigma'\in TWO(R)$.
\end{itemize}

\

\noindent Then, 
 
\begin{teo}\label{TeoAv3.5}(\cite{Avron02})\label{teoTWO} Let $\mathscr{L}$ be a sufficiently expressive language for ${\cal M}$, and let $\cal C$ be a sound and complete sequent calculus w.r.t $\cal M$. Then, $TWO({\cal C})$ is sound and complete w.r.t. ${\cal M}$ and the cut rule is admissible in $TWO({\cal C})$.
\end{teo}

\

\section{Proof system based on sets of signed formulas  for {\bf Six}}

In this section,  we shall construct both the 6-sequent calculus and its translation to a two-sided ordinary sequent calculus for \Six\, following the method described in Section \ref{s3}. As it was observed in \cite{Avron02}, the $n$-sequent calculi obtained using the above general method are hardly optimal (the same is true for the two-sided calculi). Therefore in both cases, it can be used the three general {\em streamlining principles} (which will be described in the next section) from \cite{Avron01} to reduce them to a more compact form. 

\begin{defi} The 6-sequent calculus \SFSix \, is defined as follows:\\[2mm]
{\bf Axioms}
$$ \textstyle \{\0:\alpha, \ut:\alpha, \n:\alpha, \b:\alpha, \dt:\alpha, \1:\alpha\} \mbox{ \,  where \, } \alpha \in \mathfrak{Fm}$$
{\bf Structural rules}
$$\displaystyle\frac{\Omega}{\Omega^{'}} \mbox{ \,  where \, } \Omega \subseteq\Omega^{'}$$
{\bf Logical rules} \, $i, j \in {\cal T}_6$
$$\mbox{\rm ($\vee_{i,j}$) \, } \displaystyle\frac{\Omega, i:\alpha, \, \,   j:\beta}{\Omega, \Sup \{i,j\}: \alpha \vee \beta} \hspace{2cm} \mbox{\rm ($\wedge_{i,j}$) \, } \displaystyle\frac{\Omega, i:\alpha, \, \,   j:\beta}{\Omega, \Inf \{i,j\}: \alpha \wedge \beta}$$
$$\mbox{\rm ($\neg_{i}$) \, } \displaystyle\frac{\Omega, i:\alpha}{\Omega, \1 - i: \neg\alpha} \mbox{ \, for } i\not=\n \mbox{ and }  i\not=\b $$
$$\mbox{\rm ($\neg_{\n}$) \, } \displaystyle\frac{\Omega, \n:\alpha}{\Omega, \n:\neg\alpha} \hspace{2cm} \mbox{\rm ($\neg_{\b}$) \, } \displaystyle\frac{\Omega, \b:\alpha}{\Omega, \b:\neg\alpha}$$
$$\mbox{\rm ($\nabla_{\0}$) \, } \displaystyle\frac{\Omega, \0:\alpha}{\Omega, \0:\nabla\alpha} \hspace{2cm} \mbox{\rm ($\nabla_{i}$) \, } \displaystyle\frac{\Omega, i:\alpha}{\Omega, \1:\nabla\alpha} \mbox{ \, for } i\not= \0$$
\end{defi}

\

\begin{prop}
\begin{itemize}
\item[]
\item[\rm (i)] \SFSix \ is sound and complete w.r.t the matrix  ${\cal M}_{6}$,
\item[\rm (ii)] The cut rule is admissible in \SFSix .
\end{itemize}
\end{prop}
\begin{dem} By construction and taking into account Theorem \ref{TeoAv01}. 
\end{dem}

\

\noindent Now, we translate the 6-sequent calculus \SFSix \ into a calculus of ordinary sequents. Let us call \TSSix\  the two-sided (ordinary) sequent calculus obtained after the translation depicted in Section \ref{s3}.

\begin{prop}\label{AsBs} The language of  {\bf Six} is sufficiently expressive.
\end{prop}
\begin{dem}
Indeed, it is enough to check the following.

\begin{center}
\begin{tabular}{lcl}
$v(\alpha)=\0$ & \, $\Longleftrightarrow$ \, & $v(\alpha),v(\nabla\alpha) \in {\cal N}$, $v(\neg\alpha) \in {\cal D}$ \\
$v(\alpha)=\ut$ & \, $\Longleftrightarrow$ \, &  $v(\alpha) \in {\cal N}, \,v(\neg\alpha), v(\nabla\alpha)\in {\cal D}$ \\
$v(\alpha)=\n$ & \, $\Longleftrightarrow$ \, &  $v(\alpha),v(\neg\alpha)\in {\cal N}$\\
$v(\alpha)=\b$ & \, $\Longleftrightarrow$ \, &  $v(\alpha), v(\neg\alpha)\in {\cal D}$\\
$v(\alpha)=\dt$ & \, $\Longleftrightarrow$ \,& $v(\neg\alpha) \in {\cal N}, \, v(\alpha), v(\nabla \neg\alpha)\in {\cal D}$\\
$v(\alpha)=\1$ & \, $\Longleftrightarrow$ \, & $v(\neg\alpha),  v(\nabla \neg\alpha) \in {\cal N}, \, v(\alpha)\in {\cal D}$\\
\end{tabular}
\end{center}
\end{dem}

\noindent Then, the formulas  $\alpha_j^i$'s and $\beta_k^i$'s (of Definition \ref{SufficientlyExp}) are:

\begin{center}
\begin{tabular}{|c|c|c|}  \hline
 & & \\
$i$ & $\alpha_j^i$'s & $\beta_k^i$'s \\  \hline
 & & \\
$\0$ & $\alpha_1^{\0}=p$, $\alpha_2^{\0}=\nabla p$ & $\beta_1^{\0}=\neg p$ \\ \hline
 & & \\
$\ut$ & $\alpha_1^{\ut}= p$ & $\beta_1^{\ut}= \neg p$, $\beta_2^{\ut}=\nabla p$ \\ \hline
 & & \\
$\n$ & $\alpha_1^{\n}= p$,  $\alpha_2^{\n}= \neg p$ & -- \\ \hline
 & & \\
$\b$ & -- &  $\beta_1^{\b}= p$, $\beta_2^{\b}=\neg p$\\ \hline
 & & \\
$\dt$ & $\alpha_1^{\dt}= \neg p$ & $\beta_1^{\dt}= p$, $\beta_2^{\dt}=\nabla \neg p$\\ \hline
& & \\
$\1$ & $\alpha_1^{\1}= \neg p$, $\alpha_2^{\1}=\nabla \neg p$& $\beta_1^{\1}= p$, \\ \hline
\end{tabular}
\end{center}

\

\

\noindent We start translating the axiom of \SFSix . In order to do so, we have to calculate the set of ordinary sequents
$$TWO(\alpha \, | \, \alpha \, | \,\alpha \, | \,\alpha \, | \, \alpha \, | \, \alpha)$$

\noindent This set has 324 elements and all of them can be obtained from \, $\alpha \Rightarrow \alpha$ \, using weakening\footnote{We wrote a little program in Python to calculate all of them.}.

\noindent In order to translate the logical rules, we shall take into account that

\

\

\hspace{-1.3cm}
\begin{tabular}{ll}
$TWO( \, \alpha \, | \hspace{.4cm} |\hspace{.4cm} |\hspace{.4cm}|\hspace{.4cm} |\hspace{.5cm} ) =\{ \alpha\Rightarrow , \, \nabla\alpha\Rightarrow , \, \Rightarrow\neg \alpha \} $ & $TWO( \hspace{.4cm} |\hspace{.4cm}  |\hspace{.4cm}|\, \alpha \, |\hspace{.4cm} |\hspace{.5cm} ) =\{\Rightarrow\alpha ,  \, \Rightarrow\neg\alpha\}$ \\
& \\
$TWO( \hspace{.4cm} | \, \alpha \,  |\hspace{.4cm} |\hspace{.4cm}|\hspace{.4cm} |\hspace{.5cm} ) =\{ \alpha\Rightarrow , \, \Rightarrow\neg\alpha , \, \Rightarrow\nabla\alpha \} $ & $TWO( \hspace{.4cm} |\hspace{.4cm}  |\hspace{.4cm}|\hspace{.4cm} |\, \alpha \, |\hspace{.4cm} ) =\{\neg\alpha\Rightarrow , \, \Rightarrow\alpha , \, \Rightarrow\nabla\neg\alpha\}$ \\
& \\
$TWO( \hspace{.4cm} |\hspace{.4cm}  |\, \alpha \, |\hspace{.4cm}|\hspace{.4cm} |\hspace{.5cm} ) =\{\alpha\Rightarrow ,  \, \neg\alpha\Rightarrow \} $ & $TWO( \hspace{.5cm} |\hspace{.4cm}  |\hspace{.4cm}|\hspace{.4cm} | \hspace{.4cm} |\, \alpha \, ) =\{\neg\alpha\Rightarrow ,  \, \nabla\neg\alpha\Rightarrow , \, \Rightarrow\alpha\}$ \\
\end{tabular}

\

\

\noindent Then, the ordinary-sequent rules for the negation $\neg$ are \footnote{ We have removed the contexts for simplicity.}:

\

\begin{tabular}{lll}
\\
$(\neg_{\0})_1 \,\displaystyle\frac{\alpha\Rightarrow \hspace{.3cm} \nabla\alpha\Rightarrow\hspace{.3cm} \Rightarrow\neg \alpha}{\neg\neg\alpha\Rightarrow}$ \hspace{.5cm} &
$(\neg_{\0})_2 \,  \displaystyle\frac{\alpha\Rightarrow \hspace{.3cm} \nabla\alpha\Rightarrow \hspace{.3cm} \Rightarrow\neg\alpha}{\nabla\neg\neg\alpha\Rightarrow}$ \hspace{.5cm} &
$(\neg_{\0})_3 \, \displaystyle\frac{\alpha\Rightarrow \hspace{.3cm} \nabla\alpha\Rightarrow \hspace{.3cm} \Rightarrow\neg\alpha}{\Rightarrow\neg\alpha}$\\
\\
$(\neg_{\ut})_1 \, \displaystyle\frac{\alpha\Rightarrow \hspace{.3cm} \Rightarrow\neg\alpha \hspace{.3cm} \Rightarrow\nabla\alpha}{\neg\neg\alpha\Rightarrow}$ &

$(\neg_{\ut})_2 \, \displaystyle\frac{\alpha\Rightarrow \hspace{.3cm}\Rightarrow\neg\alpha \hspace{.3cm} \Rightarrow\nabla\alpha}{\Rightarrow\neg\alpha}$ &

$(\neg_{\ut})_3 \, \displaystyle\frac{\alpha\Rightarrow \hspace{.3cm} \Rightarrow\neg\alpha \hspace{.3cm} \Rightarrow\nabla\alpha}{\Rightarrow\nabla\neg\neg\alpha}$\\
\\
\end{tabular}

$$(\neg_{\n})_1 \, \displaystyle\frac{\alpha\Rightarrow \hspace{.3cm}  \neg\alpha\Rightarrow}{\neg\alpha\Rightarrow} \hspace{1.5cm} (\neg_{\n})_2 \, \displaystyle\frac{\alpha\Rightarrow \hspace{.3cm}  \neg\alpha\Rightarrow}{\neg\neg\alpha\Rightarrow}$$

$$(\neg_{\b})_1 \, \displaystyle\frac{\Rightarrow\alpha \hspace{.3cm} \Rightarrow\neg\alpha}{\Rightarrow\neg\alpha}\hspace{1.5cm} (\neg_{\b})_2 \, \displaystyle\frac{\Rightarrow\alpha \hspace{.3cm} \Rightarrow\neg\alpha}{\Rightarrow\neg\neg\alpha}$$

\begin{tabular}{lll}
\\
$(\neg_{\dt})_1 \, \displaystyle\frac{\neg\alpha\Rightarrow \hspace{.3cm} \Rightarrow\alpha \hspace{.3cm} \Rightarrow\nabla\neg\alpha}{\neg\alpha\Rightarrow}$ & 

$(\neg_{\dt})_2 \, \displaystyle\frac{\neg\alpha\Rightarrow \hspace{.3cm} \Rightarrow\alpha \hspace{.3cm} \Rightarrow\nabla\neg\alpha}{\Rightarrow\neg\neg\alpha}$ & 

$(\neg_{\dt})_3 \, \displaystyle\frac{\neg\alpha\Rightarrow \hspace{.3cm} \Rightarrow\alpha\hspace{.3cm} \Rightarrow\nabla\neg\alpha}{\Rightarrow\nabla\neg\alpha}$\\
\\
$(\neg_{\1})_1 \, \displaystyle\frac{\neg\alpha\Rightarrow \hspace{.3cm} \nabla\neg\alpha\Rightarrow \hspace{.3cm} \Rightarrow\alpha}{\neg\alpha\Rightarrow}$ & 
$(\neg_{\1})_2 \, \displaystyle\frac{\neg\alpha\Rightarrow \hspace{.3cm} \nabla\neg\alpha\Rightarrow \hspace{.3cm} \Rightarrow\alpha}{\nabla\neg\alpha\Rightarrow}$ &
$(\neg_{\1})_3 \, \displaystyle\frac{\neg\alpha\Rightarrow \hspace{.3cm} \nabla\neg\alpha\Rightarrow \hspace{.3cm} \Rightarrow\alpha}{\Rightarrow\neg\neg\alpha}$\\
\\
\end{tabular}

\

\noindent Clearly, the rules \, $(\neg_{\0})_3$, $(\neg_{\ut})_2$, $(\neg_{\n})_1$, $(\neg_{\b})_1$,  $(\neg_{\dt})_1$ and $(\neg_{\dt})_3$,
$(\neg_{\1})_1$  and $(\neg_{\1})_2$  are superfluous (the sequent conclusion is also a sequent premise) and can be eliminated. 

\

\noindent Analogously, we obtain the rules for $\nabla$.

\

\

\begin{tabular}{lll}
$(\nabla_{\0})_1 \, \displaystyle\frac{\alpha\Rightarrow \hspace{.3cm} \nabla\alpha\Rightarrow \hspace{.3cm} \Rightarrow\neg\alpha}{\nabla\alpha\Rightarrow}$ \hspace{.5cm} &
$(\nabla_{\0})_2 \, \displaystyle\frac{\alpha\Rightarrow \hspace{.3cm} \nabla\alpha\Rightarrow  \hspace{.3cm} \Rightarrow\neg\alpha}{\nabla\nabla\alpha\Rightarrow}$ \hspace{.5cm}&
$(\nabla_{\0})_3 \, \displaystyle\frac{\alpha\Rightarrow  \hspace{.3cm} \nabla\alpha\Rightarrow  \hspace{.3cm} \Rightarrow\neg\alpha}{\Rightarrow\neg\nabla\alpha}$\\
& & \\
$(\nabla_{\ut})_1 \, \displaystyle\frac{\alpha\Rightarrow  \hspace{.3cm}\Rightarrow\neg\alpha  \hspace{.3cm} \Rightarrow\nabla\alpha}{\neg\nabla\alpha\Rightarrow}$ &
$(\nabla_{\ut})_2 \, \displaystyle\frac{\alpha\Rightarrow  \hspace{.3cm} \Rightarrow\neg\alpha  \hspace{.3cm} \Rightarrow\nabla\alpha}{\nabla\neg\nabla\alpha\Rightarrow}$ &
$(\nabla_{\ut})_3 \, \displaystyle\frac{\alpha\Rightarrow  \hspace{.3cm} \Rightarrow\neg\alpha  \hspace{.3cm} \Rightarrow\nabla\alpha}{\Rightarrow\nabla\alpha}$\\
& & \\
$(\nabla_{\n})_1 \, \displaystyle\frac{\alpha\Rightarrow  \hspace{.3cm} \neg\alpha\Rightarrow}{\neg\nabla\alpha\Rightarrow}$ &
$(\nabla_{\n})_2 \, \displaystyle\frac{\alpha\Rightarrow   \hspace{.3cm} \neg\alpha\Rightarrow}{\nabla\neg\nabla\alpha\Rightarrow}$ &
$(\nabla_{\n})_3 \, \displaystyle\frac{\alpha\Rightarrow   \hspace{.3cm} \neg\alpha\Rightarrow}{\Rightarrow\nabla\alpha}$\\
& & \\
$(\nabla_{\b})_1 \, \displaystyle\frac{\Rightarrow\alpha  \hspace{.3cm}  \Rightarrow\neg\alpha}{\neg\nabla\alpha\Rightarrow}$ &
$(\nabla_{\b})_2 \, \displaystyle\frac{\Rightarrow\alpha  \hspace{.3cm} \Rightarrow\neg\alpha}{\nabla\neg\nabla\alpha\Rightarrow}$&
$(\nabla_{\b})_3 \, \displaystyle\frac{\Rightarrow\alpha \hspace{.3cm} \Rightarrow\neg\alpha}{\Rightarrow\nabla\alpha}$\\
& & \\
$(\nabla_{\dt})_1 \, \displaystyle\frac{\neg\alpha\Rightarrow  \hspace{.3cm} \Rightarrow\alpha  \hspace{.3cm} \Rightarrow\nabla\neg\alpha}{\neg\nabla\alpha\Rightarrow}$ &
$(\nabla_{\dt})_2 \, \displaystyle\frac{\neg\alpha\Rightarrow  \hspace{.3cm} \Rightarrow\alpha  \hspace{.3cm} \Rightarrow\nabla\neg\alpha}{\nabla\neg\nabla\alpha\Rightarrow}$&
$(\nabla_{\dt})_3 \, \displaystyle\frac{\neg\alpha\Rightarrow  \hspace{.3cm} \Rightarrow\alpha  \hspace{.3cm} \Rightarrow\nabla\neg\alpha}{\Rightarrow\nabla\alpha}$\\
& & \\
$(\nabla_{\1})_1 \, \displaystyle\frac{\neg\alpha\Rightarrow  \hspace{.3cm} \nabla\neg\alpha\Rightarrow  \hspace{.3cm} \Rightarrow\alpha}{\neg\nabla\alpha\Rightarrow}$&
$(\nabla_{\1})_2 \, \displaystyle\frac{\neg\alpha\Rightarrow  \hspace{.3cm} \nabla\neg\alpha\Rightarrow  \hspace{.3cm} \Rightarrow\alpha}{\nabla\neg\nabla\alpha\Rightarrow}$&
$(\nabla_{\1})_3 \, \displaystyle\frac{\neg\alpha\Rightarrow  \hspace{.3cm} \nabla\neg\alpha\Rightarrow  \hspace{.3cm} \Rightarrow\alpha}{\Rightarrow\nabla\alpha}$\\
\end{tabular}

\

\

\noindent Again, rules $(\nabla_{\0})_1$, $(\nabla_{\0})_3$, $(\nabla_{\ut})_3$  and  $(\nabla_{\1})_3$ are superfluous and, therefore, can be eliminated.

\

\noindent The rules for  the disjunction are:

\

\

\hspace{-.9cm}
\begin{tabular}{ll}
$(\vee_{\0,\0}) \, \displaystyle\frac{\alpha\Rightarrow \hspace{.2cm} \nabla\alpha\Rightarrow \hspace{.2cm} \Rightarrow\neg\alpha \hspace{.2cm}  \beta\Rightarrow \hspace{.2cm} \nabla\beta\Rightarrow \hspace{.2cm}
 \Rightarrow \neg\beta}{\alpha\vee\beta\Rightarrow;  \, \Rightarrow\neg(\alpha\vee\beta); \, \nabla(\alpha\vee\beta)\Rightarrow }$ & $(\vee_{\0,\ut}) \,\displaystyle\frac{\alpha\Rightarrow \hspace{.2cm} \nabla\alpha\Rightarrow \hspace{.2cm} \Rightarrow\neg\alpha \hspace{.2cm}  \beta\Rightarrow\hspace{.2cm} \Rightarrow\neg\beta\hspace{.2cm} \Rightarrow \nabla\beta}{\alpha\vee\beta\Rightarrow;  \, \Rightarrow\neg(\alpha\vee\beta); \, \nabla(\alpha\vee\beta)\Rightarrow }$ \\
& \\
$(\vee_{\0,\n}) \, \displaystyle\frac{\alpha\Rightarrow \hspace{.2cm} \nabla\alpha\Rightarrow \hspace{.2cm} \Rightarrow\neg\alpha \hspace{.2cm}  \beta\Rightarrow \hspace{.2cm}  \neg\beta\Rightarrow}{(\alpha\vee\beta)\Rightarrow; \, \neg(\alpha\vee\beta)\Rightarrow}$ &
$(\vee_{\0, \b}) \displaystyle\frac{\alpha\Rightarrow \hspace{.2cm} \nabla\alpha\Rightarrow \hspace{.2cm} \Rightarrow\neg\alpha\hspace{.2cm} \Rightarrow \beta \hspace{.2cm} \Rightarrow \neg\beta}{\Rightarrow(\alpha\vee\beta); \, \Rightarrow\neg(\alpha\vee\beta)}$ \\
& \\
$(\vee_{\0,\dt}) \, \displaystyle\frac{\alpha\Rightarrow \hspace{.2cm} \nabla\alpha\Rightarrow \hspace{.2cm}  \Rightarrow\neg\alpha \hspace{.2cm} \neg\beta\Rightarrow\hspace{.2cm} \Rightarrow\beta\hspace{.2cm}\Rightarrow \nabla\neg\beta}{\Rightarrow\alpha\vee\beta; \,  \neg(\alpha\vee\beta)\Rightarrow; \, \Rightarrow\nabla\neg(\alpha\vee\beta)}$ & $(\vee_{\0,\1}) \, \displaystyle\frac{\alpha\Rightarrow \hspace{.2cm} \nabla\alpha\Rightarrow \hspace{.2cm}  \Rightarrow\neg\alpha \hspace{0,2cm} \neg\beta\Rightarrow \hspace{.2cm} \nabla\neg\beta\Rightarrow \hspace{.2cm} \Rightarrow\beta}{\Rightarrow\alpha\vee\beta; \,  \neg(\alpha\vee\beta)\Rightarrow; \, \nabla\neg(\alpha\vee\beta)\Rightarrow}$ \\
& \\
$(\vee_{\ut,\0}) \, \displaystyle\frac{\alpha\Rightarrow \hspace{.2cm} \Rightarrow\neg\alpha \hspace{.2cm} \Rightarrow \nabla\alpha\hspace{.2cm}  \beta\Rightarrow \hspace{.2cm} \nabla\beta \Rightarrow \hspace{.2cm}\Rightarrow \neg\beta}{\alpha\vee\beta\Rightarrow; \, \Rightarrow\neg(\alpha\vee\beta); \, \Rightarrow\nabla(\alpha\vee\beta)}$ & $(\vee_{\ut,\ut})\, \displaystyle\frac{ \alpha\Rightarrow \hspace{.2cm} \Rightarrow\neg\alpha \hspace{.2cm} \Rightarrow \nabla\alpha \hspace{.2cm}  \beta\Rightarrow \hspace{.2cm} \Rightarrow\neg\beta \hspace{.2cm} \Rightarrow \nabla\beta}{\alpha\vee\beta\Rightarrow; \, \Rightarrow\neg(\alpha\vee\beta); \, \Rightarrow\nabla(\alpha\vee\beta)}$ \\
&\\
$(\vee_{\ut,\n}) \, \displaystyle\frac{\alpha\Rightarrow \hspace{.2cm} \Rightarrow\neg\alpha \hspace{.2cm} \Rightarrow\nabla\alpha \hspace{.2cm}  \beta\Rightarrow \hspace{.2cm} \neg\beta\Rightarrow}{\alpha\vee\beta\Rightarrow; \, \neg(\alpha\vee\beta)\Rightarrow}$ & 
$(\vee_{\ut,\b})\displaystyle\frac{\alpha\Rightarrow \hspace{.2cm} \Rightarrow\neg\alpha \hspace{.2cm} \Rightarrow\nabla\alpha \hspace{.2cm} \Rightarrow \beta \hspace{.2cm} \Rightarrow \neg\beta}{\Rightarrow \alpha\vee\beta; \, \Rightarrow\neg(\alpha\vee\beta)}$ \\
& \\
$(\vee_{\ut,\dt}) \, \displaystyle\frac{\alpha\Rightarrow \hspace{.2cm}\Rightarrow\neg\alpha \hspace{.2cm} \Rightarrow\nabla\alpha \hspace{.2cm} \neg\beta\Rightarrow \hspace{.2cm} \Rightarrow\beta \hspace{.2cm} \Rightarrow \nabla\neg\beta}{\Rightarrow \alpha\vee\beta; \,\neg(\alpha\vee\beta)\Rightarrow; \, \Rightarrow\nabla\neg(\alpha\vee\beta)}$ & $(\vee_{\ut, \1}) \, \displaystyle\frac{\alpha\Rightarrow \hspace{.2cm} \Rightarrow\neg\alpha \hspace{.2cm} \Rightarrow\nabla\alpha \hspace{.2cm} \neg\beta\Rightarrow \hspace{.2cm} \nabla\neg\beta \Rightarrow \hspace{.2cm} \Rightarrow\beta}{\Rightarrow \alpha\vee\beta; \, \neg(\alpha\vee\beta)\Rightarrow; \,  \nabla\neg(\alpha\vee\beta)\Rightarrow}$ \\
& \\
$(\vee_{\n,\0}) \, \displaystyle\frac{\alpha\Rightarrow \hspace{.2cm} \neg\alpha\Rightarrow \hspace{.2cm}  \beta\Rightarrow \hspace{.2cm} \nabla\beta\Rightarrow \hspace{.2cm} 
 \Rightarrow \neg\beta}{\alpha\vee\beta\Rightarrow; \neg(\alpha\vee\beta)\Rightarrow}$ & $(\vee_{\n,\ut}) \, \displaystyle\frac{\alpha\Rightarrow \hspace{.2cm} \neg\alpha\Rightarrow \hspace{.2cm}  \beta\Rightarrow \hspace{.2cm} \Rightarrow\neg\beta \hspace{.2cm}  \Rightarrow \nabla\beta}{\alpha\vee\beta\Rightarrow; \, \neg(\alpha\vee\beta)\Rightarrow}$ \\
& \\
\end{tabular}
\hspace{-.9cm}
\begin{tabular}{ll}
$(\vee_{\n,\n}) \, \displaystyle\frac{\alpha\Rightarrow \hspace{.2cm}\neg\alpha\Rightarrow\hspace{.2cm} \beta\Rightarrow \hspace{.2cm} \neg\beta\Rightarrow}{\alpha\vee\beta\Rightarrow; \, \neg(\alpha\vee\beta)\Rightarrow}$ &
$(\vee_{\n,\b}) \, \displaystyle\frac{\alpha\Rightarrow \hspace{.2cm}\neg\alpha\Rightarrow\hspace{.2cm}\Rightarrow \beta \hspace{.2cm} \Rightarrow \neg\beta}{\Rightarrow \alpha\vee\beta; \, \neg(\alpha\vee\beta)\Rightarrow; \, \Rightarrow\nabla\neg(\alpha\vee\beta)}$ \\
&\\
$(\vee_{\n,\dt}) \, \displaystyle\frac{\alpha\Rightarrow \hspace{.2cm} \neg\alpha\Rightarrow\hspace{.2cm} \neg\beta\Rightarrow\hspace{.2cm}\Rightarrow\beta\hspace{.2cm}\Rightarrow \nabla\neg\beta}{\Rightarrow \alpha\vee\beta; \, \neg(\alpha\vee\beta)\Rightarrow; \, \Rightarrow\nabla\neg(\alpha\vee\beta)}$ &
$(\vee_{\n,\1}) \, \displaystyle\frac{\alpha\Rightarrow\hspace{.2cm}\neg\alpha\Rightarrow\hspace{.2cm} \neg\beta\Rightarrow \hspace{.2cm} \nabla\neg\beta \Rightarrow \beta}{\Rightarrow \alpha\vee\beta; \, \neg(\alpha\vee\beta)\Rightarrow; \, \nabla\neg(\alpha\vee\beta)\Rightarrow}$\\
&\\
$(\vee_{\b,\0}) \, \displaystyle\frac{\Rightarrow \alpha \hspace{.2cm} \Rightarrow \neg\alpha \hspace{.2cm} \beta\Rightarrow \hspace{.2cm} \nabla\beta\Rightarrow \hspace{.2cm}
 \Rightarrow \neg\beta}{\Rightarrow \alpha\vee\beta; \, \Rightarrow \neg(\alpha\vee\beta) }$ &
$(\vee_{\b, \ut}) \, \displaystyle\frac{\Rightarrow \alpha \hspace{.2cm} \Rightarrow \neg\alpha\hspace{.2cm} \beta\Rightarrow \hspace{.2cm}  \Rightarrow\neg\beta \hspace{.2cm}
 \Rightarrow \nabla\beta}{\Rightarrow \alpha\vee\beta; \, \Rightarrow \neg(\alpha\vee\beta)}$ \\
&\\
$(\vee_{\b, \n}) \, \displaystyle\frac{\Rightarrow \alpha \hspace{.2cm} \Rightarrow \neg\alpha\hspace{.2cm} \beta\Rightarrow \hspace{.2cm}  \neg\beta\Rightarrow }{\Rightarrow \alpha\vee\beta; \, \neg(\alpha\vee\beta)\Rightarrow; \, \Rightarrow\nabla\neg(\alpha\vee\beta)}$ &
$(\vee_{\b,\b})\displaystyle\frac{\Rightarrow \alpha \hspace{.2cm} \Rightarrow \neg\alpha\hspace{.2cm}\Rightarrow \beta \hspace{.2cm} \Rightarrow \neg\beta}{\Rightarrow  \alpha\vee\beta; \, \Rightarrow\neg(\alpha\vee\beta)}$ \\
& \\
$(\vee_{\b,\dt}) \, \displaystyle\frac{\Rightarrow \alpha\hspace{.2cm}\Rightarrow \neg\alpha\hspace{.2cm}\neg\beta\Rightarrow \hspace{.2cm}  \Rightarrow\beta \hspace{.2cm} \Rightarrow \nabla\neg\beta}{\Rightarrow \alpha\vee\beta; \, \neg(\alpha\vee\beta)\Rightarrow; \, \Rightarrow\nabla\neg(\alpha\vee\beta)}$ &
$(\vee_{\b,\1}) \, \displaystyle\frac{\Rightarrow \alpha \hspace{.2cm} \Rightarrow \neg\alpha\hspace{.2cm} \neg\beta\Rightarrow \hspace{.2cm} \nabla\neg\beta \Rightarrow \hspace{.2cm} \Rightarrow\beta }{\Rightarrow \alpha\vee\beta; \, \neg(\alpha\vee\beta)\Rightarrow; \, \nabla\neg(\alpha\vee\beta) \Rightarrow}$ \\
& \\
$(\vee_{\dt,\0}) \, \displaystyle\frac{ \neg\alpha\Rightarrow \hspace{.2cm}  \Rightarrow\alpha \hspace{.2cm}  \Rightarrow \nabla\neg\alpha\hspace{.2cm}\beta\Rightarrow \hspace{.2cm} \nabla\beta \Rightarrow \hspace{.2cm}  \Rightarrow \neg\beta}{\Rightarrow \alpha\vee\beta; \, \neg(\alpha\vee\beta)\Rightarrow; \, \Rightarrow \nabla\neg(\alpha\vee\beta) }$ &
$(\vee_{\dt,\ut}) \, \displaystyle\frac{\neg\alpha\Rightarrow \hspace{.2cm}\Rightarrow\alpha \hspace{.2cm} \Rightarrow\nabla\neg\alpha \hspace{.2cm} \beta\Rightarrow \hspace{.2cm} \Rightarrow \neg\beta \hspace{.2cm} \Rightarrow\nabla\beta }{\Rightarrow \alpha\vee\beta; \, \neg(\alpha\vee\beta)\Rightarrow; \, \Rightarrow \nabla\neg(\alpha\vee\beta) }$ \\
& \\
$(\vee_{\dt,\n}) \, \displaystyle\frac{\neg\alpha\Rightarrow \hspace{.2cm} \Rightarrow\alpha \hspace{.2cm} \Rightarrow\nabla\neg\alpha \hspace{.2cm}  \beta\Rightarrow \hspace{.2cm} \neg\beta\Rightarrow}{\Rightarrow  \alpha\vee\beta; \, \neg(\alpha\vee\beta)\Rightarrow; \, \Rightarrow\nabla\neg(\alpha\vee\beta)}$ & $(\vee_{\dt,\b}) \, \displaystyle\frac{\neg\alpha\Rightarrow \hspace{.2cm} \Rightarrow\alpha \hspace{.2cm} \Rightarrow\nabla\neg\alpha\hspace{.2cm}  \Rightarrow \beta \hspace{.2cm} \Rightarrow \neg\beta}{\Rightarrow  \alpha\vee\beta; \, \neg(\alpha\vee\beta)\Rightarrow; \, \Rightarrow\nabla\neg(\alpha\vee\beta)}$ 
\end{tabular}

\hspace{-1cm}
\begin{tabular}{ll}
$(\vee_{\dt,\dt}) \, \displaystyle\frac{\neg\alpha\Rightarrow \hspace{.2cm} \Rightarrow\alpha \hspace{.2cm} \Rightarrow\nabla\neg\alpha \hspace{.2cm}   \neg\beta\Rightarrow \hspace{.2cm} \Rightarrow\beta \hspace{.2cm} \Rightarrow \nabla\neg\beta}{\Rightarrow  \alpha\vee\beta; \, \neg(\alpha\vee\beta)\Rightarrow; \, \Rightarrow\nabla\neg(\alpha\vee\beta)}$ &
$(\vee_{\dt,\1}) \, \displaystyle\frac{\neg\alpha\Rightarrow \hspace{.2cm} \Rightarrow\alpha \hspace{.2cm} \Rightarrow\nabla\neg\alpha\hspace{.2cm} \neg\beta\Rightarrow \hspace{.2cm} \nabla\neg\beta\Rightarrow \hspace{.2cm} \Rightarrow\beta}{\Rightarrow  \alpha\vee\beta; \, \neg(\alpha\vee\beta)\Rightarrow; \, \nabla\neg(\alpha\vee\beta)\Rightarrow}$ \\
& \\
$(\vee_{\1,\0}) \, \displaystyle\frac{\neg\alpha\Rightarrow \hspace{.2cm} \nabla\neg\alpha\Rightarrow \hspace{.2cm}\Rightarrow\alpha\hspace{.2cm}  \beta\Rightarrow \hspace{.2cm} \nabla\beta\Rightarrow \hspace{.2cm} \Rightarrow \neg\beta}{\Rightarrow  \alpha\vee\beta; \, \neg(\alpha\vee\beta)\Rightarrow; \, \nabla\neg(\alpha\vee\beta)\Rightarrow}$ &
$(\vee_{\1, \ut}) \, \displaystyle\frac{\neg\alpha\Rightarrow \hspace{.2cm}\nabla\neg\alpha\Rightarrow \hspace{.2cm} \Rightarrow \alpha\hspace{.2cm}   \neg\beta\Rightarrow \hspace{.2cm} \nabla\neg\beta \Rightarrow \beta}{\Rightarrow  \alpha\vee\beta; \, \neg(\alpha\vee\beta)\Rightarrow; \, \nabla\neg(\alpha\vee\beta)\Rightarrow}$ \\
& \\
$(\vee_{\1, \n}) \, \displaystyle\frac{\neg\alpha\Rightarrow \hspace{.2cm} \nabla\neg\alpha\Rightarrow \hspace{.2cm} \Rightarrow \alpha\hspace{.2cm}    \beta\Rightarrow \hspace{.2cm} \neg\beta\Rightarrow }{\Rightarrow  \alpha\vee\beta; \, \neg(\alpha\vee\beta)\Rightarrow; \, \nabla\neg(\alpha\vee\beta)\Rightarrow}$ &
$(\vee_{\1, \b}) \, \displaystyle\frac{\neg\alpha\Rightarrow \hspace{.2cm} \nabla\neg\alpha\Rightarrow \hspace{.2cm}\Rightarrow\alpha\hspace{.2cm}    \Rightarrow \beta \hspace{.2cm} \Rightarrow \neg\beta}{\Rightarrow  \alpha\vee\beta; \, \neg(\alpha\vee\beta)\Rightarrow; \, \nabla\neg(\alpha\vee\beta)\Rightarrow}$ \\
& \\
$(\vee_{\1, \dt}) \, \displaystyle\frac{\neg\alpha\Rightarrow \hspace{.2cm}\nabla\neg\alpha\Rightarrow \hspace{.2cm}\Rightarrow\alpha\hspace{.2cm}   \neg\beta \Rightarrow \hspace{.2cm} \Rightarrow\beta \hspace{.2cm} \Rightarrow \nabla\neg\beta}{\Rightarrow  \alpha\vee\beta; \, \neg(\alpha\vee\beta)\Rightarrow; \, \nabla\neg(\alpha\vee\beta)\Rightarrow}$ &
$(\vee_{\1, \1}) \, \displaystyle\frac{\neg\alpha\Rightarrow \hspace{.2cm} \nabla\neg\alpha\Rightarrow \hspace{.2cm}\Rightarrow \alpha\hspace{.2cm}  \neg\beta\Rightarrow \hspace{.2cm}  \nabla\neg\beta\Rightarrow \hspace{.2cm}\Rightarrow\beta}{\Rightarrow  \alpha\vee\beta; \, \neg(\alpha\vee\beta)\Rightarrow; \, \nabla\neg(\alpha\vee\beta)\Rightarrow}$\\
\end{tabular}

\

\

\noindent Finally, the rules for  the conjunction are:

\

\hspace{-1cm}
\begin{tabular}{ll}
$(\wedge_{\0,\0}) \, \displaystyle\frac{\alpha\Rightarrow \hspace{.2cm} \nabla\alpha\Rightarrow \hspace{.2cm} \Rightarrow\neg\alpha \hspace{.2cm}  \beta\Rightarrow \hspace{.2cm} \nabla\beta\Rightarrow \hspace{.2cm}
 \Rightarrow \neg\beta}{\alpha\wedge\beta\Rightarrow;  \, \Rightarrow\neg(\alpha\wedge\beta); \, \nabla(\alpha\wedge\beta)\Rightarrow }$ & $(\wedge_{\0,\ut}) \,\displaystyle\frac{\alpha\Rightarrow \hspace{.2cm} \nabla\alpha\Rightarrow \hspace{.2cm} \Rightarrow\neg\alpha\hspace{.2cm}  \beta\Rightarrow \hspace{.2cm} \Rightarrow\neg\beta \hspace{.2cm} \Rightarrow \nabla\beta}{\alpha\wedge\beta\Rightarrow;  \, \Rightarrow\neg(\alpha\wedge\beta); \, \nabla(\alpha\wedge\beta)\Rightarrow}$ \\
& \\
$(\wedge_{\0,\n}) \, \displaystyle\frac{\alpha\Rightarrow \hspace{.2cm} \nabla\alpha\Rightarrow \hspace{.2cm} \Rightarrow\neg\alpha\hspace{.2cm}  \beta\Rightarrow \hspace{.2cm}  \neg\beta\Rightarrow}{\alpha\wedge\beta\Rightarrow;  \, \Rightarrow\neg(\alpha\wedge\beta); \, \nabla(\alpha\wedge\beta)\Rightarrow}$ &
$(\wedge_{\0, \b}) \displaystyle\frac{\alpha\Rightarrow \hspace{.2cm} \nabla\alpha\Rightarrow \hspace{.2cm} \Rightarrow\neg\alpha \hspace{.2cm} \Rightarrow \beta \hspace{.2cm} \Rightarrow \neg\beta}{\alpha\wedge\beta\Rightarrow;  \, \Rightarrow\neg(\alpha\wedge\beta); \, \nabla(\alpha\wedge\beta)\Rightarrow}$ \\
& \\
$(\wedge_{\0,\dt}) \, \displaystyle\frac{\alpha\Rightarrow \hspace{.2cm} \nabla\alpha\Rightarrow \hspace{.2cm}  \Rightarrow\neg\alpha \hspace{.2cm} \neg\beta\Rightarrow \hspace{.2cm} \Rightarrow\beta \hspace{.2cm} \Rightarrow \nabla\neg\beta}{\alpha\wedge\beta\Rightarrow;  \, \Rightarrow\neg(\alpha\wedge\beta); \, \nabla(\alpha\wedge\beta)\Rightarrow}$ & $(\wedge_{\0,\1}) \, \displaystyle\frac{\alpha\Rightarrow \hspace{.2cm} \nabla\alpha\Rightarrow \hspace{.2cm}  \Rightarrow\neg\alpha \hspace{0,2cm} \neg\beta\Rightarrow \hspace{.2cm} \nabla\neg\beta\Rightarrow \hspace{.2cm} \Rightarrow\beta}{\alpha\wedge\beta\Rightarrow;  \, \Rightarrow\neg(\alpha\wedge\beta); \, \nabla(\alpha\wedge\beta)\Rightarrow}$ \\
& \\
\end{tabular}
\hspace{-1cm}
\begin{tabular}{ll}
$(\wedge_{\ut,\0}) \, \displaystyle\frac{\alpha\Rightarrow \hspace{.2cm} \Rightarrow\neg\alpha \hspace{.2cm} \Rightarrow \nabla\alpha \hspace{.2cm}  \beta\Rightarrow \hspace{.2cm} \nabla\beta \Rightarrow \hspace{.2cm} \Rightarrow \neg\beta}{\alpha\wedge\beta\Rightarrow;  \, \Rightarrow\neg(\alpha\wedge\beta); \, \nabla(\alpha\wedge\beta)\Rightarrow}$ & $(\wedge_{\ut,\ut})\, \displaystyle\frac{ \alpha\Rightarrow \hspace{.2cm} \Rightarrow\neg\alpha \hspace{.2cm} \Rightarrow \nabla\alpha \hspace{.2cm}  \beta\Rightarrow \hspace{.2cm} \Rightarrow\neg\beta \hspace{.2cm} \Rightarrow \nabla\beta}{\alpha\wedge\beta\Rightarrow; \, \Rightarrow\neg(\alpha\wedge\beta); \, \Rightarrow\nabla(\alpha\wedge\beta)}$ \\
&\\
$(\wedge_{\ut,\n}) \, \displaystyle\frac{\alpha\Rightarrow \hspace{.2cm} \Rightarrow\neg\alpha \hspace{.2cm} \Rightarrow\nabla\alpha \hspace{.2cm}  \beta\Rightarrow \hspace{.2cm} \neg\beta\Rightarrow}{\alpha\wedge\beta\Rightarrow; \, \Rightarrow\neg(\alpha\wedge\beta); \, \Rightarrow\nabla(\alpha\wedge\beta)}$ & 
$(\wedge_{\ut,\b})\displaystyle\frac{\alpha\Rightarrow \hspace{.2cm} \Rightarrow\neg\alpha \hspace{.2cm} \Rightarrow\nabla\alpha \hspace{.2cm} \Rightarrow \beta \hspace{.2cm} \Rightarrow \neg\beta}{\alpha\wedge\beta\Rightarrow; \, \Rightarrow\neg(\alpha\wedge\beta); \, \Rightarrow\nabla(\alpha\wedge\beta)}$ \\
& \\
$(\wedge_{\ut,\dt}) \, \displaystyle\frac{\alpha\Rightarrow \hspace{.2cm} \Rightarrow\neg\alpha \hspace{.2cm} \Rightarrow\nabla\alpha \hspace{.2cm} \neg\beta\Rightarrow \hspace{.2cm} \Rightarrow\beta \hspace{.2cm} \Rightarrow \nabla\neg\beta}{\alpha\wedge\beta\Rightarrow; \, \Rightarrow\neg(\alpha\wedge\beta); \, \Rightarrow\nabla(\alpha\wedge\beta)}$ & $(\wedge_{\ut, \1}) \, \displaystyle\frac{\alpha\Rightarrow \hspace{.2cm} \Rightarrow\neg\alpha \hspace{.2cm} \Rightarrow\nabla\alpha \hspace{.2cm} \neg\beta\Rightarrow \hspace{.2cm} \nabla\neg\beta \Rightarrow \hspace{.2cm} \Rightarrow\beta}{\alpha\wedge\beta\Rightarrow; \, \Rightarrow\neg(\alpha\wedge\beta); \, \Rightarrow\nabla(\alpha\wedge\beta)}$ \\
& \\
$(\wedge_{\n,\0}) \, \displaystyle\frac{\alpha\Rightarrow \hspace{.2cm} \neg\alpha\Rightarrow \hspace{.2cm} \beta\Rightarrow \hspace{.2cm} \nabla\beta\Rightarrow \hspace{.2cm}
 \Rightarrow \neg\beta}{\alpha\wedge\beta\Rightarrow;  \, \Rightarrow\neg(\alpha\wedge\beta); \, \nabla(\alpha\wedge\beta)\Rightarrow}$ & $(\wedge_{\n,\ut}) \, \displaystyle\frac{\alpha\Rightarrow \hspace{.2cm} \neg\alpha\Rightarrow \hspace{.2cm}  \beta\Rightarrow \hspace{.2cm} \Rightarrow\neg\beta \hspace{.2cm}  \Rightarrow \nabla\beta}{\alpha\wedge\beta\Rightarrow; \, \Rightarrow\neg(\alpha\wedge\beta); \, \Rightarrow\nabla(\alpha\wedge\beta)}$ \\
& \\
$(\wedge_{\n,\n}) \, \displaystyle\frac{\alpha\Rightarrow \hspace{.2cm}\neg\alpha\Rightarrow\hspace{.2cm} \beta\Rightarrow \hspace{.2cm} \neg\beta\Rightarrow}{\alpha\wedge\beta\Rightarrow; \, \neg(\alpha\wedge\beta)\Rightarrow}$ &
$(\wedge_{\n,\b}) \, \displaystyle\frac{\alpha\Rightarrow \hspace{.2cm}\neg\alpha\Rightarrow\hspace{.2cm}\Rightarrow \beta \hspace{.2cm}\Rightarrow \neg\beta}{\alpha\wedge\beta\Rightarrow; \, \Rightarrow\neg(\alpha\wedge\beta); \, \Rightarrow\nabla(\alpha\wedge\beta)}$ \\
&\\
$(\wedge_{\n,\dt}) \, \displaystyle\frac{\alpha\Rightarrow \hspace{.2cm} \neg\alpha\Rightarrow\hspace{.2cm} \neg\beta\Rightarrow \hspace{.2cm} \Rightarrow\beta \hspace{.2cm} \Rightarrow \nabla\neg\beta}{\alpha\wedge\beta\Rightarrow; \, \neg(\alpha\wedge\beta)\Rightarrow}$ &
$(\wedge_{\n,\1}) \, \displaystyle\frac{\alpha\Rightarrow \hspace{.2cm}\neg\alpha\Rightarrow\hspace{.2cm} \neg\beta\Rightarrow \hspace{.2cm} \nabla\neg\beta \Rightarrow \hspace{.2cm}\Rightarrow \beta}{\alpha\wedge\beta\Rightarrow; \, \neg(\alpha\wedge\beta)\Rightarrow}$\\
&\\
$(\wedge_{\b,\0}) \, \displaystyle\frac{\Rightarrow \alpha \hspace{.2cm} \Rightarrow \neg\alpha \hspace{.2cm}  \beta\Rightarrow \hspace{.2cm} \nabla\beta\Rightarrow \hspace{.2cm}
 \Rightarrow \neg\beta}{\alpha\wedge\beta\Rightarrow;  \, \Rightarrow\neg(\alpha\wedge\beta); \, \nabla(\alpha\wedge\beta)\Rightarrow}$ &
$(\wedge_{\b, \ut}) \, \displaystyle\frac{\Rightarrow \alpha \hspace{.2cm} \Rightarrow \neg\alpha\hspace{.2cm}  \beta\Rightarrow \hspace{.2cm}  \Rightarrow\neg\beta \hspace{.2cm}
 \Rightarrow \nabla\beta}{\alpha\wedge\beta\Rightarrow; \, \Rightarrow\neg(\alpha\wedge\beta); \, \Rightarrow\nabla(\alpha\wedge\beta)}$ \\
\end{tabular}

\hspace{-1.5cm}
\begin{tabular}{ll}
$(\wedge_{\b, \n}) \, \displaystyle\frac{\Rightarrow \alpha \hspace{.2cm} \Rightarrow \neg\alpha\hspace{.2cm}\beta\Rightarrow \hspace{.2cm}  \neg\beta\Rightarrow }{\alpha\wedge\beta\Rightarrow; \, \Rightarrow\neg(\alpha\wedge\beta); \, \Rightarrow\nabla(\alpha\wedge\beta)}$ &
$(\wedge_{\b,\b})\displaystyle\frac{\Rightarrow \alpha \hspace{.2cm} \Rightarrow \neg\alpha\hspace{.2cm}\Rightarrow \beta \hspace{.2cm} \Rightarrow \neg\beta}{\Rightarrow  \alpha\wedge\beta; \, \Rightarrow\neg(\alpha\wedge\beta)}$ \\
& \\
$(\wedge_{\b,\dt}) \, \displaystyle\frac{\Rightarrow \alpha\hspace{.2cm}\Rightarrow \neg\alpha\hspace{.2cm}\neg\beta\Rightarrow \hspace{.2cm}  \Rightarrow\beta \hspace{.2cm} \Rightarrow \nabla\neg\beta}{\Rightarrow  \alpha\wedge\beta; \, \Rightarrow\neg(\alpha\wedge\beta)}$ &
$(\wedge_{\b,\1}) \, \displaystyle\frac{\Rightarrow \alpha \hspace{.2cm} \Rightarrow \neg\alpha\hspace{.2cm} \neg\beta\Rightarrow \hspace{.2cm} \nabla\neg\beta \Rightarrow \hspace{.2cm} \Rightarrow\beta }{\Rightarrow  \alpha\wedge\beta; \, \Rightarrow\neg(\alpha\wedge\beta)}$ \\
& \\
$(\wedge_{\dt,\0}) \, \displaystyle\frac{ \neg\alpha\Rightarrow \hspace{.2cm}  \Rightarrow\alpha \hspace{.2cm}  \Rightarrow \nabla\neg\alpha\hspace{.2cm}  \beta\Rightarrow \hspace{.2cm} \nabla\beta \Rightarrow \hspace{.2cm}  \Rightarrow \neg\beta}{\alpha\wedge\beta\Rightarrow;  \, \Rightarrow\neg(\alpha\wedge\beta); \, \nabla(\alpha\wedge\beta)\Rightarrow }$ &
$(\wedge_{\dt,\ut}) \, \displaystyle\frac{\neg\alpha\Rightarrow \hspace{.2cm}\Rightarrow\alpha \hspace{.2cm} \Rightarrow\nabla\neg\alpha \hspace{.2cm} \beta\Rightarrow \hspace{.2cm} \Rightarrow \neg\beta \hspace{.2cm} \Rightarrow\nabla\beta }{\alpha\wedge\beta\Rightarrow; \, \Rightarrow\neg(\alpha\wedge\beta); \, \Rightarrow\nabla(\alpha\wedge\beta) }$ \\
& \\
$(\wedge_{\dt,\n}) \, \displaystyle\frac{\neg\alpha\Rightarrow \hspace{.2cm} \Rightarrow\alpha \hspace{.2cm} \Rightarrow\nabla\neg\alpha \hspace{.2cm}  \beta\Rightarrow \hspace{.2cm} \neg\beta\Rightarrow}{\alpha\wedge\beta\Rightarrow; \, \neg(\alpha\wedge\beta)\Rightarrow}$ & $(\wedge_{\dt,\b}) \, \displaystyle\frac{\neg\alpha\Rightarrow \hspace{.2cm} \Rightarrow\alpha \hspace{.2cm} \Rightarrow\nabla\neg\alpha\hspace{.2cm}  \Rightarrow \beta \hspace{.2cm} \Rightarrow \neg\beta}{\Rightarrow  \alpha\wedge\beta; \, \Rightarrow\neg(\alpha\wedge\beta)}$ \\
& \\
$(\wedge_{\dt,\dt}) \, \displaystyle\frac{\neg\alpha\Rightarrow \hspace{.2cm} \Rightarrow\alpha \hspace{.2cm} \Rightarrow\nabla\neg\alpha \hspace{.2cm}   \neg\beta\Rightarrow \hspace{.2cm} \Rightarrow\beta \hspace{.2cm} \Rightarrow \nabla\neg\beta}{\Rightarrow  \alpha\wedge\beta; \, \neg(\alpha\wedge\beta)\Rightarrow; \, \Rightarrow\nabla\neg(\alpha\wedge\beta)}$ &
$(\wedge_{\dt,\1}) \, \displaystyle\frac{\neg\alpha\Rightarrow \hspace{.2cm} \Rightarrow\alpha \hspace{.2cm} \Rightarrow\nabla\neg\alpha\hspace{.2cm} \neg\beta\Rightarrow \hspace{.2cm} \nabla\neg\beta\Rightarrow \hspace{.2cm} \Rightarrow\beta}{\Rightarrow  \alpha\wedge\beta; \, \neg(\alpha\wedge\beta)\Rightarrow; \, \Rightarrow\nabla\neg(\alpha\wedge\beta)}$ \\
& \\
$(\wedge_{\1,\0}) \, \displaystyle\frac{\neg\alpha\Rightarrow \hspace{.2cm} \nabla\neg\alpha\Rightarrow \hspace{.2cm}\Rightarrow\alpha\hspace{.2cm}   \beta\Rightarrow \hspace{.2cm} \nabla\beta\Rightarrow \hspace{.2cm} \Rightarrow \neg\beta}{\alpha\wedge\beta\Rightarrow;  \, \Rightarrow\neg(\alpha\wedge\beta); \, \nabla(\alpha\wedge\beta)\Rightarrow }$ &
$(\wedge_{\1, \ut}) \, \displaystyle\frac{\neg\alpha\Rightarrow \hspace{.2cm}\nabla\neg\alpha\Rightarrow \hspace{.2cm} \Rightarrow \alpha\hspace{.2cm}   \beta\Rightarrow \hspace{.2cm}  \Rightarrow \neg \beta \hspace{.2cm} \Rightarrow \nabla \beta}{\alpha\wedge\beta\Rightarrow; \, \Rightarrow\neg(\alpha\wedge\beta); \, \Rightarrow\nabla(\alpha\wedge\beta)}$ \\
& \\
$(\wedge_{\1, \n}) \, \displaystyle\frac{\neg\alpha\Rightarrow \hspace{.2cm} \nabla\neg\alpha\Rightarrow \hspace{.2cm} \Rightarrow \alpha\hspace{.2cm}    \beta\Rightarrow \hspace{.2cm} \neg\beta\Rightarrow }{\alpha\wedge\beta\Rightarrow; \, \neg(\alpha\wedge\beta)\Rightarrow}$ &
$(\wedge_{\1, \b}) \, \displaystyle\frac{\neg\alpha\Rightarrow \hspace{.2cm} \nabla\neg\alpha\Rightarrow \hspace{.2cm}\Rightarrow\alpha\hspace{.2cm}     \Rightarrow \beta\hspace{.2cm}\Rightarrow \neg\beta}{\Rightarrow  \alpha\wedge\beta; \, \Rightarrow\neg(\alpha\wedge\beta)}$ \\
& \\
$(\wedge_{\1, \dt}) \, \displaystyle\frac{\neg\alpha\Rightarrow \hspace{.2cm}\nabla\neg\alpha\Rightarrow \hspace{.2cm}\Rightarrow\alpha\hspace{.2cm}  \Rightarrow  \neg\beta\Rightarrow \hspace{.2cm} \Rightarrow\beta \hspace{.2cm} \Rightarrow \nabla\neg\beta}{\Rightarrow  \alpha\wedge\beta; \, \neg(\alpha\wedge\beta)\Rightarrow; \, \Rightarrow\nabla\neg(\alpha\wedge\beta)}$ &
$(\wedge_{\1, \1}) \, \displaystyle\frac{\neg\alpha\Rightarrow \hspace{.2cm} \nabla\neg\alpha\Rightarrow \hspace{.2cm}\Rightarrow \alpha\hspace{.2cm}  \neg\beta\Rightarrow\hspace{.2cm} \nabla\neg\beta\Rightarrow \hspace{.2cm} \Rightarrow\beta}{\Rightarrow  \alpha\wedge\beta; \, \neg(\alpha\wedge\beta)\Rightarrow; \, \nabla\neg(\alpha\wedge\beta)\Rightarrow}$\\
\end{tabular}

\

\begin{prop}\label{corcompTSSix}
\begin{itemize}
\item[]
\item[\rm (i)] \TSSix \ is sound and complete w.r.t the matrix  ${\cal M}_{6}$,
\item[\rm (ii)] The cut rule is admissible in \TSSix .
\end{itemize}
\end{prop}
\begin{dem} By construction of the system \TSSix\ and Theorem \ref{TeoAv3.5}.
\end{dem}

\

\section{Streamlining the two-sided sequential system \TSSix\ }

The two-sided calculus \TSSix\  obtained using the general method of \cite{Avron01} (and reviewed in Section \ref{s3}) is as a rule hardly optimal (as it is usually the case with this type of “generic” systems).  Therefore, we use the three general streamlining principles established in \cite{Avron01} to reduce them to a more compact form. Of these three, the first and the third decrease the number of rules (which is our main measure of complexity), while the second simplifies a rule by decreasing the number of its premises (since the third rule increases this number, its application is often followed by applications of the first two). \\[2mm]
Recall that a rule (r) is context-free if whenever \, $\displaystyle\frac{\Gamma_1\Rightarrow\Delta_1 \, \dots \, \Gamma_k\Rightarrow\Delta_k}{\Gamma\Rightarrow\Delta}$ \,  is a valid application of (r), and $\Gamma'$ and $\Delta'$ are sets of formulas, then 
$$\displaystyle\frac{\Gamma_1,\Gamma'\Rightarrow\Delta_1, \Delta' \, \dots \, \Gamma_k, \Gamma' \Rightarrow\Delta_k, \Delta'}{\Gamma, \Gamma' \Rightarrow\Delta, \Delta'}$$
 is also a valid application of (r) \footnote{ Rule ($\neg$) of system $\mathfrak{S}$ in page 5 is not context-free.}. Then, the three streamlining principles are: 

\

\noindent {\bf Principle 1.} If a rule in $R$ is derivable from other rules, it can be deleted.

\

\noindent {\bf Principle 2.} If $\frac{S}{\Sigma}$ (where S is a set of premises) is a rule in $R$, $S'$ is a subset of $S$ and $\frac{S'}{\Sigma}$ is derivable in $R$ (perhaps using cuts), then $\frac{S}{\Sigma}$ can be replaced by $\frac{S'}{\Sigma}$. In particular: if $\frac{S}{\Sigma}$  is a rule in $R$, $\Gamma\Rightarrow \Delta \in S$, and $\Gamma\Rightarrow \Delta$ is derivable from $S\setminus \{\Gamma\Rightarrow \Delta\}$ in $R$, then $\frac{S}{\Sigma}$ can be replaced with $\frac{S\setminus \{\Gamma\Rightarrow \Delta\}}{\Sigma}$. Two very simple, but
quite useful cases of this are when $\Gamma\Rightarrow \Delta$ is subsumed by an axiom or by some sequent in $\Gamma\Rightarrow \Delta$.

\

\noindent {\bf Principle 3.} Let $G$ be a sequent system. If in $G$ we have the rules $\displaystyle \frac{\{\Gamma_{i}\Rightarrow\Delta_{i}\}_{\small 1 \leq i \leq n}} {\Gamma\Rightarrow\Delta}$ , $\displaystyle \frac{\{\Pi_{j}\Rightarrow\Theta_{j}\}_{1 \leq j\leq k} }{\Gamma\Rightarrow\Delta}$ and both are context--free, then we can replace these two rules by the new rule   
$$\displaystyle \frac{\{\Gamma_{i},\Pi_{j}\Rightarrow\Delta_{i},\Theta_{j}\}_{\small  1\leq i \leq n,  1 \leq j \leq k }} {\Gamma\Rightarrow\Delta}$$

\

The proofs of Principles 1 and 2 are immediate. Principle 3 is obtained using weakening appropriately, the contraction rule (which is implicit since we are dealing with set of formulas instead of sequences of formulas) and the fact that both rules are context free.
Along with these principles, we shall use systematically the following proposition.

\

\begin{prop}\label{propred} Let $\mathfrak{SC}$ be a sequent calculus having the identity axiom and the structural rules of weakening and contraction (left and right); and suppose that $\displaystyle \frac{S \cup \{\Rightarrow  \varphi\}}{\Sigma}$ and $\displaystyle \frac{S \cup \{\varphi\Rightarrow\}}{\Sigma}$ are two context-free rules of $\mathfrak{SC}$. Then,  $\displaystyle \frac{S}{\Sigma}$ is derivable in $\mathfrak{SC}$.
\end{prop}
\begin{dem} Suppose that $S=\{\Gamma_1\Rightarrow \Delta_1, \dots, \Gamma_n\Rightarrow \Delta_n\}$ and $\Sigma$ is $\Gamma\Rightarrow \Delta$. Since ${\rm (r1)} \, \displaystyle \frac{S \cup \{\Rightarrow  \varphi\}}{\Sigma}$ and ${\rm (r2)}\ , \displaystyle \frac{S \cup \{\varphi\Rightarrow\}}{\Sigma}$ are context-free rules, by Principle 3, we have that

$${\rm (r)} \, \displaystyle \frac{\{\Gamma_{i},\Gamma_{j}\Rightarrow\Delta_{i},\Delta_{j}\}_{\small  1\leq i,j \leq n} \cup \{\varphi, \Gamma_{i}\Rightarrow\Delta_{i}\}_{\small  1\leq i \leq n} \cup \{ \Gamma_{i}\Rightarrow\Delta_{i}, \varphi\}_{\small  1\leq i \leq n} \cup \{\varphi\Rightarrow\varphi \}} {\Gamma\Rightarrow\Delta}$$
is derivable in $\mathfrak{SC}$. Now, for every $i\in\{1, \dots, n\}$, the sequents $\Gamma_{i},\Gamma_{j}\Rightarrow\Delta_{i},\Delta_{j}$, $1\leq j \leq n$; \, $\varphi, \Gamma_i \Rightarrow \Delta_i$ and $\Gamma_i \Rightarrow \Delta_i, \varphi$ can be obtained from $\Gamma_i \Rightarrow \Delta_i$ via the weakening rules, which is an upper sequent of (r) (for $i=j$). Besides, $\varphi\Rightarrow\varphi$ is an instance of the identity axiom. Therefore,   $\displaystyle \frac{S}{\Sigma}$ is derivable in $\mathfrak{SC}$.
\end{dem}

\

\noindent From rules $(\neg_{\0})_1$ and $(\neg_{\ut})_1$ and Proposition \ref{propred}, we obtain the following derivable rule

$$ \displaystyle\frac{\alpha\Rightarrow \, \, \Rightarrow\neg \alpha}{\neg\neg\alpha\Rightarrow}$$

\noindent and, from this rule, the rule $(\neg_{\n})_2$ and Proposition \ref{propred}, we have that

$$\displaystyle\frac{\alpha\Rightarrow }{\neg\neg\alpha\Rightarrow}$$

\noindent is derivable in \SFSix . Now, after restoring context, we have

$$\mbox{($\neg\neg\Rightarrow$)} \, \, \displaystyle\frac{ \alpha, \Gamma \Rightarrow \Delta}{\neg\neg\alpha, \Gamma\Rightarrow \Delta}$$

\noindent In an analogous way, from rules $(\neg_{\dt})_2$,  $(\neg_{\1})_3$, $(\neg_{\b})_2$ and Proposition  \ref{propred} we get

$$\mbox{($\Rightarrow\neg\neg$)} \, \, \displaystyle\frac{\Gamma \Rightarrow \Delta,  \alpha}{ \Gamma\Rightarrow \Delta, \neg\neg\alpha}$$
Therefore, the rules for $\neg$ are
$$\mbox{($\neg\neg\Rightarrow$)} \, \, \displaystyle\frac{ \alpha, \Gamma \Rightarrow \Delta}{\neg\neg\alpha, \Gamma\Rightarrow \Delta} \hspace{1.5cm} \mbox{($\Rightarrow\neg\neg$)} \, \, \displaystyle\frac{\Gamma \Rightarrow \Delta,  \alpha}{ \Gamma\Rightarrow \Delta, \neg\neg\alpha}$$
$$(\nabla\neg\neg\Rightarrow)  \, \displaystyle\frac{\alpha, \Gamma\Rightarrow \Delta \hspace{.5cm} \nabla\alpha, \Gamma\Rightarrow \Delta  \hspace{.5cm} \Gamma \Rightarrow\Delta,\neg\alpha}{\nabla\neg\neg\alpha, \Gamma\Rightarrow \Delta}\hspace{1cm}(\Rightarrow\nabla\neg\neg)  \, \displaystyle\frac{\alpha, \Gamma\Rightarrow \Delta \hspace{.5cm}  \Gamma\Rightarrow \Delta, \neg\alpha  \hspace{.5cm} \Gamma \Rightarrow\Delta,\nabla\alpha}{\Gamma\Rightarrow \Delta, \nabla\neg\neg\alpha}$$

\

\noindent Concerning the connective $\nabla$, we have that from rules  $(\nabla_{\dt})_3$,  $(\nabla_{\1})_3$  and Proposition \ref{propred} we get $\, \displaystyle\frac{\Rightarrow \alpha , \, \neg\alpha\Rightarrow}{\Rightarrow\nabla\alpha}$. From this rule, using $(\nabla_{\b})_3$ and Proposition \ref{propred} we get \, $\displaystyle\frac{\Rightarrow\alpha}{\Rightarrow\nabla\alpha}$; and using $(\nabla_{\n})_3$ and Proposition \ref{propred}, we obtain  \,  $\displaystyle\frac{\neg\alpha\Rightarrow}{\Rightarrow\nabla\alpha}$.\\
On the other hand, from $(\nabla_{\1})_2$, $(\nabla_{\dt})_2$ and Proposition \ref{propred} we obtain $(1) \displaystyle\frac{\neg\alpha\Rightarrow  \,  \Rightarrow\alpha}{\nabla\neg\nabla\alpha\Rightarrow}$. We combine this rule in two ways: with $(\nabla_{\n})_2$ and Proposition \ref{propred} we have  $ (2)\displaystyle\frac{\neg\alpha\Rightarrow }{\nabla\neg\nabla\alpha\Rightarrow}$; and with  $(\nabla_{\b})_2$ and Proposition \ref{propred} we have $\displaystyle\frac{\Rightarrow\alpha}{\nabla\neg\nabla\alpha\Rightarrow}$. Now, we combine rules (1), (2) and $(\nabla_{\ut})_2$ using  Principle 3 and obtain 

$$(\nabla\neg\nabla\Rightarrow) \, \displaystyle\frac{\Gamma, \neg\alpha\Rightarrow  \Delta, \alpha, \nabla\alpha}{\Gamma, \nabla\neg\nabla\alpha\Rightarrow \Delta}$$
Reasoning similarly and combining the appropriate rules we have the following derivable rules (after restoring contex) in \SFSix .

$$ (\Rightarrow\nabla) \, \displaystyle\frac{\Gamma, \neg \alpha\Rightarrow\Delta, \alpha}{\Gamma\Rightarrow\Delta, \nabla\alpha}  \hspace{2cm} (\nabla\nabla\Rightarrow) \,  \displaystyle\frac{\alpha, \Gamma\Rightarrow \Delta \hspace{.5cm} \nabla\alpha, \Gamma\Rightarrow\Delta \hspace{.5cm} \Gamma \Rightarrow\Delta, \neg\alpha}{\nabla\nabla\alpha, \Gamma\Rightarrow\Delta}$$
$$ (\neg\nabla\Rightarrow) \, \displaystyle\frac{\Gamma, \neg\alpha\Rightarrow  \Delta, \alpha, \nabla\alpha}{\neg\nabla\alpha,\Gamma\Rightarrow\Delta} \hspace{2cm} (\Rightarrow\neg\nabla) \, \displaystyle\frac{\alpha, \Gamma\Rightarrow \Delta  \hspace{.5cm} \nabla\alpha, \Gamma\Rightarrow\Delta \hspace{.5cm} \Gamma\Rightarrow\Delta, \neg\alpha}{\Gamma\Rightarrow\neg\nabla\alpha, \Delta}$$
$$ (\nabla\neg\nabla\Rightarrow) \, \displaystyle\frac{\Gamma, \neg\alpha\Rightarrow  \Delta, \alpha, \nabla\alpha}{\Gamma, \nabla\neg\nabla\alpha\Rightarrow \Delta}$$

\

\noindent On the other hand, we have the following:

\begin{prop} The rules \, $\displaystyle\frac{\nabla \alpha\Rightarrow}{ \alpha\Rightarrow}$,  \, $\displaystyle\frac{\nabla \alpha\Rightarrow}{\Rightarrow\neg\alpha}$ \, and \, $\displaystyle\frac{\Rightarrow\nabla \alpha}{\Rightarrow\nabla \neg\neg\alpha}$ \, are admissible in \TSSix\ .
\end{prop}
\begin{dem}  We shall prove it only for the first rule (the others are analogous). Suppose that the sequent $\nabla \alpha\Rightarrow$ is provable in \TSSix . Then, by Proposition  \ref{corcompTSSix}, $\models_{{\cal M}_6} \nabla \alpha\Rightarrow$. Then, for every ${\cal M}_6$-valuation $v$, $v(\nabla \alpha)=\hat{\nabla} v(\alpha)\in \{\0,\ut,\n\}$. By the definition of $\hat{\nabla}$ (see Section 2), we have that $\hat{\nabla} v(\alpha)=\0$; and then, $v(\alpha)=\0$, for every valuation $v$. That is, $\models_{{\cal M}_6} \alpha\Rightarrow$; and by Proposition  \ref{corcompTSSix}, $\alpha\Rightarrow$ is provable in \TSSix .
\end{dem}

\

\noindent Then, from this and  Principles 2 and 3, we have that the rules ($\nabla\neg\neg \Rightarrow$), ($\Rightarrow\nabla\neg\neg$), ($\nabla\nabla\Rightarrow$) and ($\Rightarrow\neg\nabla$) can be recast as follows:

$$(\nabla\neg\neg\Rightarrow)  \, \displaystyle\frac{\nabla\alpha, \Gamma\Rightarrow \Delta}{\nabla\neg\neg\alpha, \Gamma\Rightarrow \Delta}\hspace{2cm}(\Rightarrow\nabla\neg\neg)  \, \displaystyle\frac{\Gamma \Rightarrow\Delta,\nabla\alpha}{\Gamma\Rightarrow \Delta, \nabla\neg\neg\alpha}$$

$$(\nabla\nabla\Rightarrow) \,  \displaystyle\frac{\nabla\alpha, \Gamma\Rightarrow\Delta}{\nabla\nabla\alpha, \Gamma\Rightarrow\Delta} \hspace{2cm} (\Rightarrow\neg\nabla) \, \displaystyle\frac{\nabla\alpha, \Gamma\Rightarrow\Delta}{\Gamma\Rightarrow\neg\nabla\alpha, \Delta}$$

\

\noindent Reasoning analogously and combining (using exclusively Proposition \ref{propred}) the appropriate rules we can obtain the following derivable rules in  \SFSix\ for the connective $\vee$ .\\[2mm]
From $(\vee_{\0,\0})$ and $(\vee_{\0,\ut})$  we get \,  $\displaystyle\frac{\alpha\Rightarrow \, \nabla\alpha\Rightarrow \, \Rightarrow\neg\alpha \hspace{0,5cm}  \beta\Rightarrow \, \Rightarrow \neg\beta}{(\alpha\vee\beta)\Rightarrow}$. From this and $(\vee_{\0,\n})$ we get \,  $\mbox{(I) } \, \displaystyle\frac{\alpha\Rightarrow \, \nabla\alpha\Rightarrow \, \Rightarrow\neg\alpha \hspace{0,5cm}  \beta\Rightarrow}{(\alpha\vee\beta)\Rightarrow}$. Besides,  from $(\vee_{\ut,\0})$ and $(\vee_{\ut,\ut})$ we get \,  $\displaystyle\frac{\alpha\Rightarrow \,\Rightarrow\nabla\alpha \, \Rightarrow\neg\alpha \hspace{0,5cm} \beta\Rightarrow   \Rightarrow \neg\beta}{(\alpha\vee\beta)\Rightarrow}$; \, and from this and $(\vee_{\ut,\n})$ we get \, $\displaystyle\frac{\alpha\Rightarrow \,\Rightarrow\nabla\alpha \, \Rightarrow\neg\alpha \hspace{0,5cm} \beta\Rightarrow}{(\alpha\vee\beta)\Rightarrow}$. From this last rule and (I) we obtain \,   (II) \, $\displaystyle\frac{\alpha\Rightarrow \, \Rightarrow\neg\alpha \hspace{0,5cm} \beta\Rightarrow}{(\alpha\vee\beta)\Rightarrow}$. \\
On the other hand, from $(\vee_{\n, \0})$ and $(\vee_{\n, \ut})$ we have  \,  $\displaystyle\frac{\alpha\Rightarrow \, \neg\alpha\Rightarrow \hspace{0,5cm} \beta\Rightarrow \, \Rightarrow\neg\beta}{(\alpha\vee\beta)\Rightarrow}$ \, and combining it with  $(\vee_{\n,\n})$ we get \, $\displaystyle\frac{\alpha\Rightarrow,\neg\alpha\Rightarrow \hspace{0,5cm}  \beta\Rightarrow}{(\alpha\vee\beta)\Rightarrow}$. From this last rule and (II) we obtain \,  $\displaystyle\frac{\alpha\Rightarrow \hspace{0,5cm}  \beta\Rightarrow}{(\alpha\vee\beta)\Rightarrow}$. Then, after restoring context, we have the following (classical) derivable rule in \TSSix .
$$\mbox{ ($\vee\Rightarrow$) }\displaystyle\frac{\alpha, \Gamma\Rightarrow\Delta \hspace{0,5cm}  \beta, \Gamma\Rightarrow\Delta}{\alpha\vee\beta
, \Gamma\Rightarrow\Delta}$$

\noindent Similarly, from $(\vee_{\0,\0})$ and $(\vee_{\0,\ut})$: \, $\displaystyle\frac{\alpha\Rightarrow \,\nabla\alpha\Rightarrow \, \Rightarrow\neg\alpha \hspace{0,5cm}  \beta\Rightarrow \, \Rightarrow \neg\beta}{\Rightarrow\neg(\alpha\vee\beta)}$. From this and $(\vee_{\0,\b})$: \, $ \mbox{(III)} \, \displaystyle\frac{\alpha\Rightarrow \, \nabla\alpha\Rightarrow \, \Rightarrow\neg\alpha \hspace{0,5cm}  \Rightarrow \neg\beta}{\Rightarrow\neg(\alpha\vee\beta)}$. On the other hand, from $(\vee_{\ut,\0})$ and $(\vee_{\ut,\ut})$ we obtain \\ $\displaystyle\frac{\alpha\Rightarrow \, \Rightarrow\neg\alpha \, \Rightarrow\nabla\alpha \hspace{0,5cm}  \beta\Rightarrow \, \Rightarrow \neg\beta}{\Rightarrow\neg(\alpha\vee\beta)}$; and combining this one with $(\vee_{\ut,\b})$: $ \displaystyle\frac{\alpha\Rightarrow \, \Rightarrow\neg\alpha \, \Rightarrow\nabla\alpha \hspace{0,5cm}  \Rightarrow \neg\beta}{\Rightarrow\neg(\alpha\vee\beta)}$. From  this one and (III) we get \, $\mbox{(IV)} \,\displaystyle\frac{\alpha\Rightarrow \,\Rightarrow\neg\alpha \hspace{0,5cm} \Rightarrow \neg\beta}{\Rightarrow\neg(\alpha\vee\beta)}$. On the other hand,  from
$(\vee_{\b,\0})$ and $(\vee_{\b,\ut})$ we get \, $\displaystyle\frac{\Rightarrow\alpha, \Rightarrow\neg\alpha \hspace{0,5cm}  \beta\Rightarrow, \Rightarrow \neg\beta}{\Rightarrow\neg(\alpha\vee\beta)}$; from this and $(\vee_{\b,\b})$ we have \, $\displaystyle\frac{\Rightarrow\alpha, \Rightarrow\neg\alpha \hspace{0,5cm}   \Rightarrow \neg\beta}{\Rightarrow\neg(\alpha\vee\beta)}$. Combining this last rule with (IV) we get \,  $\displaystyle\frac{\Rightarrow\neg\alpha \hspace{0,5cm}   \Rightarrow \neg\beta}{\Rightarrow\neg(\alpha\vee\beta)}$. Now, restoring context we obtain the next derivable rule in \TSSix .
$$\mbox{ ($\Rightarrow\neg\vee$) } \, \,\displaystyle\frac{\Gamma\Rightarrow\neg\alpha, \Delta \hspace{0,5cm}   \Gamma\Rightarrow \neg\beta, \Delta}{\Gamma\Rightarrow\neg(\alpha\vee\beta), \Delta}$$
 
\noindent This last rule is highly expected to be derivable in \TSSix \ since involutive Stone algebras are, in particular, De Morgan algebras and therefor the De Morgan laws hold. \\
The rule $\displaystyle\frac{\Rightarrow\beta }{\Rightarrow \alpha\vee\beta}$ needs a little more work to be deduced. We display its deduction process in the Table \ref{tabla1}. \\[2mm]
Similarly, we can obtain  \, $\displaystyle\frac{\Rightarrow\alpha }{\Rightarrow \alpha\vee\beta}$. By Principle 3 and restoring context we have
$$\mbox{ ($\Rightarrow\vee$) } \, \, \displaystyle\frac{\Gamma\Rightarrow\alpha, \beta, \Delta }{\Gamma\Rightarrow \alpha\vee\beta, \Delta}$$
Following an analogous reasoning we can obtain the rule
$$\mbox{ ($\neg\vee\Rightarrow$) } \, \, \displaystyle\frac{\neg\alpha, \neg\beta, \Gamma\Rightarrow \Delta }{\neg(\alpha\vee\beta), \Gamma\Rightarrow \Delta}$$

\begin{table}
\begin{center}
\begin{tabular}{|c|c|c|}  \hline
Step & Rules to be combined & Resulting rule \\ \hline \hline 
& & \\[-1.5mm]
(1) &$(\vee_{\1,\dt})$ , \,  $(\vee_{\1,\1})$ & $\displaystyle\frac{\neg\alpha\Rightarrow \, \nabla\neg\alpha\Rightarrow \, \Rightarrow \alpha\hspace{0,5cm}      \Rightarrow\beta \, \neg\beta \Rightarrow }{\Rightarrow \alpha\vee\beta}$ \\ \hline
& &  \\[-1.5mm]
(2) & $(\vee_{\dt,\dt})$ , \,  $(\vee_{\dt,\1})$ & $\displaystyle\frac{\Rightarrow\alpha \, \neg\alpha\Rightarrow \, \Rightarrow\nabla\neg\alpha\hspace{0,5cm}      \Rightarrow\beta \, \neg\beta \Rightarrow }{\Rightarrow \alpha\vee\beta}$\\ \hline
& &  \\[-1.5mm]
(3) & $(1)$ , \,  $(2)$ & $\displaystyle\frac{\Rightarrow\alpha \,\neg\alpha\Rightarrow\hspace{0,5cm}      \Rightarrow\beta \, \neg\beta \Rightarrow }{\Rightarrow \alpha\vee\beta}$\\ \hline
& & \\
(4) & $(\vee_{\1,\n})$, \,  $(\vee_{\dt,\n})$ & $\displaystyle\frac{\Rightarrow\alpha \, \neg\alpha\Rightarrow\hspace{0,5cm} \beta\Rightarrow \, \neg\beta \Rightarrow}{\Rightarrow \alpha\vee\beta}$\\ \hline
& &  \\[-1.5mm]
(5) & $(\vee_{\1,\b})$, \,  $(\vee_{\dt,\b})$ & $\displaystyle\frac{\Rightarrow\alpha \, \neg\alpha\Rightarrow\hspace{0,5cm}      \Rightarrow\beta \, \Rightarrow \neg\beta }{\Rightarrow \alpha\vee\beta}$\\ \hline
& &  \\[-1.5mm]
(6) & $(3)$, \,  $(5)$ & $\displaystyle\frac{\Rightarrow\alpha \, \neg\alpha\Rightarrow\hspace{0,5cm}      \Rightarrow\beta }{\Rightarrow \alpha\vee\beta}$ \\ \hline
& &  \\[-1.5mm]
(7) & $(\vee_{\n,\1})$, \,   $(\vee_{\n,\dt})$ & $\displaystyle\frac{\alpha\Rightarrow \, \neg\alpha\Rightarrow\hspace{0,5cm} \Rightarrow\beta \, \neg\beta \Rightarrow}{\Rightarrow \alpha\vee\beta}$\\  \hline
& &  \\[-1.5mm]
(8) & $(7)$, \,   $(\vee_{\n,\b})$ & $\displaystyle\frac{\alpha\Rightarrow \,\neg\alpha\Rightarrow\hspace{0,5cm} \Rightarrow\beta}
{\Rightarrow \alpha\vee\beta}$ \\  \hline
& &  \\[-1.5mm]
(9) & $(6)$, \,   $(8)$ & $\displaystyle\frac{\neg\alpha\Rightarrow\hspace{0,5cm} \Rightarrow\beta}
{\Rightarrow \alpha\vee\beta}$ \\  \hline
& &  \\[-1.5mm]
(10) & $(\vee_{\b,\1})$, \, $(\vee_{\b,\dt})$& $\displaystyle\frac{\Rightarrow\alpha \, \Rightarrow\neg\alpha\hspace{0,5cm}      \Rightarrow\beta \, \neg\beta\Rightarrow }{\Rightarrow \alpha\vee\beta}$ \\  \hline
& &  \\[-1.5mm]
(11) & $(10)$, \, $(\vee_{\b,\b})$& $\displaystyle\frac{\Rightarrow\alpha \, \Rightarrow\neg\alpha \hspace{0,5cm}      \Rightarrow\beta }{\Rightarrow \alpha\vee\beta}$ \\  \hline
& &  \\[-1.5mm]
(12) & $(\vee_{\0,\1})$, \, $(\vee_{\0,\dt})$ &$\displaystyle\frac{\alpha\Rightarrow \, \Rightarrow\neg\alpha \, \nabla\alpha\Rightarrow \hspace{0,5cm}      \Rightarrow\beta \, \neg\beta\Rightarrow}{\Rightarrow \alpha\vee\beta}$ \\  \hline
& & \\[-1.5mm]
(13) & $(\vee_{\ut,\1})$, \,  $(\vee_{\ut,\dt})$ &$\displaystyle\frac{\alpha\Rightarrow \, \Rightarrow\neg\alpha \,\Rightarrow\nabla\alpha\hspace{0,5cm}\Rightarrow\beta \, \neg\beta \Rightarrow }{\Rightarrow \alpha\vee\beta}$ \\  \hline
& &  \\[-1.5mm]
(14) & $(12)$, \,  $(13)$ &$\displaystyle\frac{\alpha\Rightarrow \, \Rightarrow\neg\alpha\hspace{0,5cm}\Rightarrow\beta \, \neg\beta \Rightarrow }{\Rightarrow \alpha\vee\beta}$ \\  \hline
& & \\
(15) & $(\vee_{\0,\b})$, \,  $(\vee_{\ut,\b})$ &$\displaystyle\frac{\alpha\Rightarrow \, \Rightarrow\neg\alpha\hspace{0,5cm}      \Rightarrow\beta \, \Rightarrow\neg\beta }{\Rightarrow \alpha\vee\beta}$ \\  \hline
& & \\[-1.5mm]
(16) & $(14)$, \,  $(15)$ &$\displaystyle\frac{\alpha\Rightarrow \,\Rightarrow\neg\alpha\hspace{0,5cm}\Rightarrow\beta }{\Rightarrow \alpha\vee\beta}$ \\  \hline
& & \\[-1.5mm]
(17) & $(11)$, \,  $(16)$ &$\displaystyle\frac{\Rightarrow\neg\alpha \hspace{0,5cm}\Rightarrow\beta }{\Rightarrow \alpha\vee\beta}$ \\  \hline
& & \\[-1.5mm]
(18) & $(9)$, \,  $(17)$ &$\displaystyle\frac{\Rightarrow\beta }{\Rightarrow \alpha\vee\beta}$\\  \hline
\end{tabular}
\end{center}
\caption{Derivation of  $\displaystyle\frac{\Rightarrow\beta }{\Rightarrow \alpha\vee\beta}$}
\label{tabla1}
\end{table}

\begin{prop} The following rules are admissible in \TSSix .
\begin{center}
\begin{tabular}{ll}
{\rm (i)} \, $\displaystyle \frac{\nabla \alpha \Rightarrow \hspace{0.5cm } \nabla \beta \Rightarrow }{\nabla (\alpha\vee\beta) \Rightarrow}$ \hspace{2cm} & {\rm (ii)} \, $\displaystyle \frac{\Rightarrow\nabla \alpha,  \nabla \beta }{\Rightarrow\nabla (\alpha\vee\beta)}$ \\
& \\
{\rm (iii)} \, $\displaystyle \frac{\nabla \neg\alpha, \nabla\neg\beta \Rightarrow }{\nabla \neg (\alpha\vee\beta) \Rightarrow}$ \hspace{2cm} & {\rm (iv)} \, $\displaystyle \frac{\Rightarrow\nabla \neg \alpha \hspace{0.5cm }  \Rightarrow\nabla\neg \beta }{\Rightarrow\nabla \neg(\alpha\vee\beta)}$ \\
\end{tabular}
\end{center}
\end{prop}
\begin{dem} (i): Suppose that $\nabla \alpha \Rightarrow$ and $\nabla \beta \Rightarrow$ are provable. By Proposition \ref{corcompTSSix}, $\models_{{\cal M}_6} \nabla \alpha \Rightarrow$ and $\models_{{\cal M}_6} \nabla \beta \Rightarrow$. Let $v$ be a valuation. Then, $v(\nabla \alpha)=\hat{\nabla}v(\alpha)=\0$ (since $\hat{\nabla}x = \0$ iff $x=\0$ for every $x\in {\cal T}_6$) and $v( \nabla \beta)=\hat{\nabla}v(\beta)=\0$. Hence, $v(\nabla(\alpha\vee\beta))=\hat{\nabla}(v(\alpha)\hat{\vee}v(\beta))= \hat{\nabla}v(\alpha)\hat{\vee}\hat{\nabla}v(\beta)=\0\vee\0=\0$ (recall that $\nabla$ is a lattice-homomorphism in the class \Sto ).\\[2mm]
(ii): Similar to (i) and taking into account that in ${\cal M}_6$ it holds: $x\not=\0$ or $y\not=\0$ imply $x \hat{\vee} y\not=\0$.\\[2mm]
(iii): Similar to (i) and taking into account that in ${\cal M}_6$ it holds: $\hat{\neg}x=\0$ or $\hat{\neg}y=\0$ imply $\hat{\neg}(x \hat{\vee} y)=\0$.\\[2mm]
(iv):  Similar to (i) and taking into account that in ${\cal M}_6$ it holds: $\hat{\neg}x\not=\0$ and $\hat{\neg}y\not=\0$ imply $\hat{\neg}(x \hat{\vee} y)\not=\0$.
\end{dem}

\

\noindent Therefore,  the streamlined rules for the disjunction are:
$$\mbox{ ($\vee\Rightarrow$) }\displaystyle\frac{\alpha, \Gamma\Rightarrow\Delta \hspace{0,5cm}  \beta, \Gamma\Rightarrow\Delta}{\alpha\vee\beta
, \Gamma\Rightarrow\Delta} \hspace{2cm} \mbox{ ($\Rightarrow\vee$) } \, \, \displaystyle\frac{\Gamma\Rightarrow\alpha, \beta, \Delta }{\Gamma\Rightarrow \alpha\vee\beta, \Delta}$$
$$\mbox{ ($\neg\vee\Rightarrow$) } \, \, \displaystyle\frac{\neg\alpha, \neg\beta, \Gamma\Rightarrow \Delta }{\neg(\alpha\vee\beta), \Gamma\Rightarrow \Delta} \hspace{2cm} \mbox{ ($\Rightarrow\neg\vee$) } \, \,\displaystyle\frac{\Gamma\Rightarrow\neg\alpha, \Delta \hspace{0,5cm}   \Gamma\Rightarrow \neg\beta, \Delta}{\Gamma\Rightarrow\neg(\alpha\vee\beta), \Delta} $$
$$\mbox{ ($\nabla\vee\Rightarrow$) } \, \displaystyle \frac{\nabla \alpha, \Gamma \Rightarrow \Delta \hspace{0.5cm } \nabla \beta, \Gamma \Rightarrow \Delta}{\nabla (\alpha\vee\beta), \Gamma \Rightarrow \Delta} \hspace{2cm} \mbox{ ($\Rightarrow\nabla\vee$) }  \, \frac{\Gamma \Rightarrow\nabla \alpha,  \nabla \beta, \Delta}{\Gamma\Rightarrow\nabla (\alpha\vee\beta), \Delta} $$
$$\mbox{ ($\nabla\neg\vee\Rightarrow$) } \, \displaystyle \frac{\nabla \neg\alpha, \nabla\neg\beta, \Gamma \Rightarrow \Delta }{\nabla \neg (\alpha\vee\beta), \Gamma \Rightarrow\Delta} \hspace{2cm} \mbox{ ($\Rightarrow\nabla\neg\vee$) }  \,  \frac{\Gamma\Rightarrow\nabla \neg \alpha, \Delta \hspace{0.5cm }  \Gamma\Rightarrow\nabla\neg \beta, \Delta }{\Gamma\Rightarrow\nabla \neg(\alpha\vee\beta), \Delta} $$

\

\noindent The streamlining process for the rules of $\wedge$ goes similarly to the above; and the resulting rules are the next.
$$\mbox{ ($\wedge\Rightarrow$) }\displaystyle\frac{\alpha, \beta,\Gamma\Rightarrow\Delta}{\alpha\wedge\beta
, \Gamma\Rightarrow\Delta} \hspace{2cm} \mbox{ ($\Rightarrow\wedge$) } \, \, \displaystyle\frac{\Gamma\Rightarrow\alpha, \Delta  \hspace{0,5cm} \Gamma\Rightarrow\beta, \Delta  }{\Gamma\Rightarrow \alpha\wedge\beta, \Delta}$$
$$\mbox{ ($\neg\wedge\Rightarrow$) } \, \, \displaystyle\frac{\neg\alpha,  \Gamma\Rightarrow \Delta \hspace{0,5cm} \neg\beta, \Gamma\Rightarrow \Delta}{\neg(\alpha\wedge\beta), \Gamma\Rightarrow \Delta} \hspace{2cm} \mbox{ ($\Rightarrow\neg\wedge$) } \, \,\displaystyle\frac{\Gamma\Rightarrow\neg\alpha, \neg\beta, \Delta}{\Gamma\Rightarrow\neg(\alpha\wedge\beta), \Delta} $$
$$\mbox{ ($\nabla\wedge\Rightarrow$) } \, \displaystyle \frac{\nabla \alpha, \nabla \beta,\Gamma \Rightarrow \Delta}{\nabla (\alpha\wedge\beta), \Gamma \Rightarrow \Delta} \hspace{2cm} \mbox{ ($\Rightarrow\nabla\wedge$) }  \, \frac{\Gamma \Rightarrow\nabla \alpha,  \Delta \hspace{0,5cm} \Gamma \Rightarrow  \nabla \beta, \Delta }{\Gamma\Rightarrow\nabla (\alpha\wedge\beta), \Delta} $$
$$\mbox{ ($\nabla\neg\wedge\Rightarrow$) } \, \displaystyle \frac{\nabla \neg\alpha,  \Gamma \Rightarrow \Delta \hspace{0,5cm}  \nabla\neg\beta, \Gamma \Rightarrow \Delta}{\nabla \neg (\alpha\wedge\beta), \Gamma \Rightarrow\Delta} \hspace{2cm} \mbox{ ($\Rightarrow\nabla\neg\wedge$) }  \,  \frac{\Gamma\Rightarrow\nabla \neg \alpha, \nabla\neg \beta, \Delta }{\Gamma\Rightarrow\nabla \neg(\alpha\wedge\beta), \Delta} $$

\

\noindent Finally, we propose a new principle, Principle 4, which will allow us to further optimize some rules.

\begin{prop}(Principle 4)  Let $\mathfrak{SC}$ be a sequent calculus having the structural rules of weakening and contraction (left and right). If ${\rm (r)} \, \displaystyle \frac{\{\Gamma_{i}\Rightarrow\Delta_{i}\}_{\small  1\leq i \leq n}} {\Gamma\Rightarrow\Delta}$ is a rule of $\mathfrak{SC}$ and there are sets of formulas $\Gamma'_i$'s and $\Delta'_i$ for $1\leq i\leq n$ such that:
\begin{itemize}
\item[(i)] $\Gamma'_i\subseteq \Gamma_i$ and $\Delta'_i\subseteq \Delta'$,
\item[(ii)] the sequent $\Gamma'_i\Rightarrow\Delta'_i$ is derivable in $\mathfrak{SC}$ from $\Gamma_i\Rightarrow \Delta_i$ without using the (r) rule nor the cut rule, from $\Gamma_i\Rightarrow \Delta_i$,
\end{itemize}
Then, (r) can be replaced by \, ${\rm (r')} \, \displaystyle \frac{\{\Gamma'_{i}\Rightarrow\Delta'_{i}\}_{\small  1\leq i \leq n}} {\Gamma\Rightarrow\Delta}$.
\end{prop} 
\begin{dem}
Let us see that {\rm (r)} and {\rm (r')} are equivalent. It is clear that from (r) we can prove (r') using suitable weakenings. Conversely, from (r') we can prove (r). Indeed, 
\begin{prooftree}
\AxiomC{$\Gamma_1 \Rightarrow \Delta_1$}
\noLine
\UnaryInfC{$\vdots$}
\noLine
\UnaryInfC{$\Gamma'_1 \Rightarrow \Delta'_1$}

\AxiomC{$\dots$}
\noLine
\UnaryInfC{}
\noLine
\UnaryInfC{}

\AxiomC{$\Gamma_n \Rightarrow \Delta_n$}
\noLine
\UnaryInfC{$\vdots$}
\noLine
\UnaryInfC{$\Gamma'_n \Rightarrow \Delta'_n$}

\LeftLabel{\small(r')}
\TrinaryInfC{$\Gamma\Rightarrow \Delta$}
\end{prooftree}

\end{dem}

\

\noindent It is clear that if $\mathfrak{SC}$ is a cut--free sequent calculus, then the resulting system by applying Principle 4 is also cut--free.\\ 
Then, by Principle 4, we can replace the rules: 
$$\displaystyle\frac{\Gamma, \neg\alpha\Rightarrow  \Delta, \alpha, \nabla\alpha}{\neg\nabla\alpha,\Gamma\Rightarrow\Delta}\hspace{1cm} \mbox{ \, \, by \, \,}\hspace{1cm} \displaystyle\frac{\Gamma, \Rightarrow  \Delta, \nabla\alpha}{\neg\nabla\alpha,\Gamma\Rightarrow\Delta}$$

\begin{table}[H]
{\bf Axioms}
$$ \alpha \Rightarrow\alpha$$

{\bf Structural Rules}

$$ (w\Rightarrow) \, \displaystyle\frac{\Gamma\Rightarrow\Delta}{\Gamma, \alpha\Rightarrow\Delta} \hspace{2cm} (\Rightarrow w) \, \displaystyle\frac{\Gamma\Rightarrow\Delta}{\Gamma\Rightarrow\Delta, \alpha} $$

{\bf Logic Rules}

$$\mbox{ ($\vee\Rightarrow$) }\displaystyle\frac{\alpha, \Gamma\Rightarrow\Delta \hspace{0,5cm}  \beta, \Gamma\Rightarrow\Delta}{\alpha\vee\beta
, \Gamma\Rightarrow\Delta} \hspace{2cm} \mbox{ ($\Rightarrow\vee$) } \, \, \displaystyle\frac{\Gamma\Rightarrow\alpha, \beta, \Delta }{\Gamma\Rightarrow \alpha\vee\beta, \Delta}$$
$$\mbox{ ($\neg\vee\Rightarrow$) } \, \, \displaystyle\frac{\neg\alpha, \neg\beta, \Gamma\Rightarrow \Delta }{\neg(\alpha\vee\beta), \Gamma\Rightarrow \Delta} \hspace{2cm} \mbox{ ($\Rightarrow\neg\vee$) } \, \,\displaystyle\frac{\Gamma\Rightarrow\neg\alpha, \Delta \hspace{0,5cm}   \Gamma\Rightarrow \neg\beta, \Delta}{\Gamma\Rightarrow\neg(\alpha\vee\beta), \Delta} $$
$$\mbox{ ($\nabla\vee\Rightarrow$) } \, \displaystyle \frac{\nabla \alpha, \Gamma \Rightarrow \Delta \hspace{0.5cm } \nabla \beta, \Gamma \Rightarrow \Delta}{\nabla (\alpha\vee\beta), \Gamma \Rightarrow \Delta} \hspace{2cm} \mbox{ ($\Rightarrow\nabla\vee$) }  \, \frac{\Gamma \Rightarrow\nabla \alpha,  \nabla \beta, \Delta}{\Gamma\Rightarrow\nabla (\alpha\vee\beta), \Delta} $$
$$\mbox{ ($\nabla\neg\vee\Rightarrow$) } \, \displaystyle \frac{\nabla \neg\alpha, \nabla\neg\beta, \Gamma \Rightarrow \Delta }{\nabla \neg (\alpha\vee\beta), \Gamma \Rightarrow\Delta} \hspace{2cm} \mbox{ ($\Rightarrow\nabla\neg\vee$) }  \,  \frac{\Gamma\Rightarrow\nabla \neg \alpha, \Delta \hspace{0.5cm }  \Gamma\Rightarrow\nabla\neg \beta, \Delta }{\Gamma\Rightarrow\nabla \neg(\alpha\vee\beta), \Delta} $$

$$\mbox{ ($\wedge\Rightarrow$) }\displaystyle\frac{\alpha, \beta,\Gamma\Rightarrow\Delta}{\alpha\wedge\beta
, \Gamma\Rightarrow\Delta} \hspace{2cm} \mbox{ ($\Rightarrow\wedge$) } \, \, \displaystyle\frac{\Gamma\Rightarrow\alpha, \Delta  \hspace{0,5cm} \Gamma\Rightarrow\beta, \Delta  }{\Gamma\Rightarrow \alpha\wedge\beta, \Delta}$$
$$\mbox{ ($\neg\wedge\Rightarrow$) } \, \, \displaystyle\frac{\neg\alpha,  \Gamma\Rightarrow \Delta \hspace{0,5cm} \neg\beta, \Gamma\Rightarrow \Delta}{\neg(\alpha\wedge\beta), \Gamma\Rightarrow \Delta} \hspace{2cm} \mbox{ ($\Rightarrow\neg\wedge$) } \, \,\displaystyle\frac{\Gamma\Rightarrow\neg\alpha, \neg\beta, \Delta}{\Gamma\Rightarrow\neg(\alpha\wedge\beta), \Delta} $$
$$\mbox{ ($\nabla\wedge\Rightarrow$) } \, \displaystyle \frac{\nabla \alpha, \nabla \beta,\Gamma \Rightarrow \Delta}{\nabla (\alpha\wedge\beta), \Gamma \Rightarrow \Delta} \hspace{2cm} \mbox{ ($\Rightarrow\nabla\wedge$) }  \, \frac{\Gamma \Rightarrow\nabla \alpha,  \Delta \hspace{0,5cm} \Gamma \Rightarrow  \nabla \beta, \Delta }{\Gamma\Rightarrow\nabla (\alpha\wedge\beta), \Delta} $$
$$\mbox{ ($\nabla\neg\wedge\Rightarrow$) } \, \displaystyle \frac{\nabla \neg\alpha,  \Gamma \Rightarrow \Delta \hspace{0,5cm}  \nabla\neg\beta, \Gamma \Rightarrow \Delta}{\nabla \neg (\alpha\wedge\beta), \Gamma \Rightarrow\Delta} \hspace{2cm} \mbox{ ($\Rightarrow\nabla\neg\wedge$) }  \,  \frac{\Gamma\Rightarrow\nabla \neg \alpha, \nabla\neg \beta, \Delta }{\Gamma\Rightarrow\nabla \neg(\alpha\wedge\beta), \Delta} $$

$$\mbox{($\neg\neg\Rightarrow$)} \, \, \displaystyle\frac{ \alpha, \Gamma \Rightarrow \Delta}{\neg\neg\alpha, \Gamma\Rightarrow \Delta} \hspace{1.5cm} \mbox{($\Rightarrow\neg\neg$)} \, \, \displaystyle\frac{\Gamma \Rightarrow \Delta,  \alpha}{ \Gamma\Rightarrow \Delta, \neg\neg\alpha}$$
$$(\nabla\neg\neg\Rightarrow)  \, \displaystyle\frac{\nabla\alpha, \Gamma\Rightarrow \Delta}{\nabla\neg\neg\alpha, \Gamma\Rightarrow \Delta}\hspace{2cm}(\Rightarrow\nabla\neg\neg)  \, \displaystyle\frac{\Gamma \Rightarrow\Delta,\nabla\alpha}{\Gamma\Rightarrow \Delta, \nabla\neg\neg\alpha}$$

$$ (\Rightarrow\nabla) \, \displaystyle\frac{\Gamma, \neg\alpha\Rightarrow\Delta, \alpha}{\Gamma\Rightarrow\Delta, \nabla\alpha}  \hspace{2cm} (\nabla\Rightarrow) \,  \displaystyle\frac{\nabla\alpha, \Gamma\Rightarrow\Delta}{\nabla\nabla\alpha, \Gamma\Rightarrow\Delta}$$
$$ (\neg\nabla\Rightarrow) \, \displaystyle\frac{\Gamma\Rightarrow  \Delta, \nabla\alpha}{\neg\nabla\alpha,\Gamma\Rightarrow\Delta} \hspace{2cm} (\Rightarrow\neg\nabla) \, \displaystyle\frac{\nabla\alpha, \Gamma\Rightarrow\Delta}{\Gamma\Rightarrow\neg\nabla\alpha, \Delta}$$
$$ (\nabla\neg\nabla\Rightarrow) \, \displaystyle\frac{\Gamma\Rightarrow  \Delta, \nabla\alpha}{\Gamma, \nabla\neg\nabla\alpha\Rightarrow \Delta}$$

\caption{The sequent calculus \GSix\ }
\label{figgsix}
\end{table}

and
$$\displaystyle\frac{\Gamma, \neg\alpha\Rightarrow  \Delta, \alpha, \nabla\alpha}{\Gamma, \nabla\neg\nabla\alpha\Rightarrow \Delta}\hspace{1cm} \mbox{ \, \, by \, \,}\hspace{1cm} \displaystyle\frac{\Gamma\Rightarrow  \Delta, \nabla\alpha}{\Gamma, \nabla\neg\nabla\alpha\Rightarrow \Delta}$$

\

\noindent In this way, we obtain the streamlined calculus, which we name  \GSix, whose rules are displayed in Table \ref{figgsix}.  

\noindent The above construction and the detailed description of the streamlining process constitute a proof for the following theorem.

\begin{teo}\label{teocutfree}\begin{itemize}
\item[]
\item[\rm (i)] \GSix \ is sound and complete w.r.t the matrix  ${\cal M}_{6}$,
\item[\rm (ii)] The cut rule is admissible in \GSix .
\end{itemize}
\end{teo}

\

\noindent The following is a cut-free \GSix\ -proof of the sequent \, $p\vee q \Rightarrow \neg(\neg p \vee \neg q)$ (see Section \ref{s2}).

\begin{prooftree}
\AxiomC{$p \Rightarrow p$}
\LeftLabel{\small$\Rightarrow$(w)}
\UnaryInfC{$p \Rightarrow p, \neg \neg q$}
\LeftLabel{\small($\Rightarrow\neg\neg$)}
\UnaryInfC{$p\Rightarrow  \neg \neg p, \neg \neg q$}

\AxiomC{$q \Rightarrow q$}
\LeftLabel{\small$\Rightarrow$(w)}
\UnaryInfC{$q \Rightarrow  \neg \neg p, q$}
\LeftLabel{\small($\Rightarrow\neg\neg$)}
\UnaryInfC{$q\Rightarrow  \neg \neg p, \neg \neg q$}

\LeftLabel{\small($\vee\Rightarrow$)}
\BinaryInfC{$p\vee q\Rightarrow  \neg \neg p, \neg \neg q$}
\LeftLabel{\small($\Rightarrow\neg\wedge$)}
\UnaryInfC{$p\vee q\Rightarrow  \neg( \neg p \wedge \neg q)$}
\end{prooftree}

\

We end this section showing an application of the cut-elimination property. More precisely, we shall prove that, despite the fact that \Six\ is a paraconsistent logic, it entails no contradictions. This will be derived in a very simple way from the cut-elimination theorem. First, consider the following lemma.

\begin{lem}\label{cont} Let $\gamma$ be a formula. The following conditions are equivalent.
\begin{enumerate}[\rm (i)]
\item the sequent \, $\Rightarrow \gamma \wedge \neg \gamma$ \,  is provable in \GSix .
\item the empty sequent \, $\Rightarrow$ \, is provable in \GSix .
\end{enumerate}
\end{lem}
\begin{dem} (ii) $\Longrightarrow$ (i): Immediate.\\
(i) $\Longrightarrow$ (ii): Suppose that $\Rightarrow \gamma \wedge \neg \gamma$  is provable in \GSix . Then, by Theorem \ref{teocutfree}, there exist a cut-free derivation $\cal D$ of $\Rightarrow \gamma \wedge \neg \gamma$. Then, $\cal D$ is of the form
\begin{prooftree}
\AxiomC{$\vdots$}
\UnaryInfC{$\Rightarrow \gamma$}

\AxiomC{$\vdots$}
\LeftLabel{\small(r)}
\UnaryInfC{$\Rightarrow \neg \gamma$}

\LeftLabel{\small($\Rightarrow\wedge$)}
\BinaryInfC{$\Rightarrow \gamma \wedge\neg \gamma$}
\end{prooftree}
We use induction on the complexity, $c(\gamma)$, of $\gamma$.\\
{\bf Base step.} If $c(\gamma)=0$, then $\gamma$ is a propositional variable. Then, the lower sequent of inference (r) \, is $\Rightarrow \neg p$. Analyzing the rules of \GSix\  (except the cut rule), we conclude that (r) has to be right weakening and the upper sequent of this inference has to be the empty sequent.\\
{\bf Inductive step.} Suppose that $c(\gamma)>0$. Then, (r) can be one of the rules ($\Rightarrow\neg\vee$), ($\Rightarrow\neg\wedge$), ($\Rightarrow\neg\neg$) and ($\Rightarrow\neg\nabla$). We shall just analyze the case when (r) is ($\Rightarrow\neg\vee$), the others are analogous. In this case, $\cal D$ is of the following form
\begin{prooftree}
\AxiomC{(1)}
\noLine
\UnaryInfC{$\vdots$}
\UnaryInfC{$\Rightarrow \alpha, \beta$}
\LeftLabel{\small($\rightarrow\vee$)}
\UnaryInfC{$\Rightarrow \alpha\vee \beta$}

\AxiomC{(2)}
\noLine
\UnaryInfC{$\vdots$}
\UnaryInfC{$\Rightarrow  \neg\alpha$}

\AxiomC{(3)}
\noLine
\UnaryInfC{$\vdots$}
\UnaryInfC{$\Rightarrow  \neg\beta$}
\LeftLabel{\small($\rightarrow\neg\vee$)}
\BinaryInfC{$\Rightarrow  \neg(\alpha\vee \beta)$}

\LeftLabel{\small($\Rightarrow\wedge$)}
\BinaryInfC{$\Rightarrow (\alpha\vee \beta)\wedge \neg(\alpha\vee \beta)$}
\end{prooftree}
Then, we can construct the following derivation
\begin{prooftree}
\AxiomC{(1)}
\noLine
\UnaryInfC{$\vdots$}
\UnaryInfC{$\Rightarrow \alpha, \beta$}
\LeftLabel{\small($\Rightarrow$w)}
\UnaryInfC{$\Rightarrow \alpha, \beta, \neg\alpha$}

\AxiomC{(2)}
\noLine
\UnaryInfC{$\vdots$}
\UnaryInfC{$\Rightarrow  \neg\alpha$}
\LeftLabel{\small(w)'s}
\UnaryInfC{$\beta \Rightarrow \alpha, \neg\alpha$}

\LeftLabel{\small(cut)}
\BinaryInfC{$\Rightarrow \alpha, \neg\alpha$}
\LeftLabel{\small($\Rightarrow\wedge$)}
\UnaryInfC{$\Rightarrow \alpha\wedge \neg\alpha$}
\end{prooftree}
That is,  \, $\Rightarrow \alpha\wedge \neg\alpha$ is provable in \GSix . Since $c(\alpha)<c(\gamma)$, by the inductive hypothesis we have that  the empty sequent $\Rightarrow$ is provable.
\end{dem}

\

\begin{prop}\label{empty} The empty sequent is not provable in \GSix .
\end{prop}
\begin{dem} From Theorem \ref{teocutfree} and analyzing the rules of \GSix .
\end{dem}

\

\begin{cor} \Six\ does not entail contradictions.
\end{cor}
\begin{dem} From Lema \ref{cont}, Proposition \ref{empty} and Theorem \ref{teocutfree}.
\end{dem}

\section{Decision procedure à la Gentzen for \GSix}

The results presented in the previous sections were communicated in the {\em 1st Meeting Brazil-Colombia in Logic}. There,  we were asked about a decision procedure for \Six\ (based on the system \GSix ) since, as it is well known, the fact that a logic has a cut--free sequent version, per se,  does not guarantee the existence of a decision procedure for it. In this section,  we state a decision procedure for \GSix\, written in pseudocode, which is an adaptation of the one proposed by G. Gentzen for the propositional fragment of $\bf LK$ (see \cite[Theorem 6.13]{Tau}).

\begin{defi}\label{defgensub} ({\it Generalized  subformula}) The set of generalized subformulas of a given formula $\gamma$, $\gsub{\gamma}$, is defined as the least set of formulas fulfilling the following conditions:
\begin{enumerate}[(1)]
\item $\alpha \in \gsub{\alpha}$;
\item $\gsub{\alpha}\cup\gsub{\beta}\subseteq \gsub{\alpha \ast \beta} $, $\ast \in \{\vee, \wedge \}$;
\item $\gsub{\neg \alpha}\cup\gsub{\neg \beta}\subseteq \gsub{\neg (\alpha \ast \beta)}$, $\ast \in \{\vee, \wedge \}$; 
\item $\gsub{\alpha}\cup\gsub{\neg \alpha}\subseteq\gsub{\nabla \alpha}$;  
\item $\gsub{\nabla\alpha}\cup\gsub{\nabla \beta}\subseteq \gsub{\nabla(\alpha\wedge\beta)}$;  
\item $\gsub{\nabla\neg\alpha}\cup\gsub{\nabla\neg\beta}\subseteq \gsub{\nabla\neg(\alpha\wedge\beta)}$;
\item $\gsub{\alpha}\subseteq \gsub{\neg\neg \alpha}$;
\item $\gsub{\nabla\alpha}\subseteq \gsub{\nabla\neg\neg \alpha}$;
\item $\gsub{\nabla\alpha}\subseteq \gsub{\neg\nabla \alpha}$;
\end{enumerate} 
\end{defi}

\begin{rem}\label{regsub} Note that, from (4) and (9)  of Definition \ref{defgensub} we have that, for a given $\alpha$, $\alpha$ and $\neg \alpha$ are generalized subformulas of $\neg\nabla \alpha$ and
$$\gsub{\alpha}\cup\gsub{\neg\alpha}\cup\gsub{\nabla\alpha}\subseteq\gsub{\nabla\neg \nabla\alpha},$$
$$\gsub{\nabla\alpha}\subseteq\gsub{\nabla\nabla \alpha},$$
$$\gsub{\alpha}\cup\gsub{\neg\alpha}\cup\gsub{\nabla\alpha}\subseteq \gsub{\neg\nabla \alpha}.$$
\end{rem}

\

\noindent If $\Gamma=\{\gamma_1, \dots, \gamma_n\}$ ($n<\omega$), the set of all generalized subformulas of $\Gamma$ is defined as $\gsub{\Gamma}=\bigcup \limits_{i=1}^{n} \gsub{\gamma_i}$. Besides, the set of all generalized subformulas of the sequent $\Gamma\Rightarrow\Delta$ is given by 
$$\gsub{\Gamma\Rightarrow\Delta}\stackrel{def}{=}\gsub{\Gamma}\cup\gsub{\Delta}.$$
Then, we have the next

\begin{prop}  (Generalized Subformula Property) Let ${\cal D}$ be a cut--free deduction in \GSix\ of a sequent $\Gamma\Rightarrow\Delta$.
Then, for any sequent $\Gamma'\Rightarrow\Delta'$ in ${\cal D}$ we have
$$\gsub{\Gamma'\Rightarrow\Delta'} \subseteq \gsub{\Gamma\Rightarrow\Delta}.$$
\end{prop}
\begin{dem} By induction on the number inferences in ${\cal D}$, Definition \ref{defgensub}, Remark \ref{regsub} and inspection of the rules of \GSix.
\end{dem}

\begin{eje} The sequent $\nabla\neg\nabla(p\wedge q)\Rightarrow\neg\nabla p \vee \neg \nabla q$ is provable in \GSix. Then,\\[1mm]
$\gsub{\nabla\neg\nabla(p\wedge q)\Rightarrow\neg\nabla p \vee \neg \nabla}=\gsub{\nabla\neg\nabla(p\wedge q)}\cup\gsub{\neg\nabla p \vee \neg \nabla }=\{\nabla\neg\nabla(p\wedge q), \break \neg\nabla(p\wedge q),  \neg\neg\nabla(p\wedge q),  \neg\nabla p \vee \neg \nabla q, \nabla(p\wedge q), \neg(p\wedge q), p\wedge q, \neg\nabla p, \neg \nabla q, \nabla p, \nabla q,  \neg p, \neg q, p, q \}$; and  there exists a cut--free deduction of it. Indeed

\begin{prooftree}
\AxiomC{$\nabla p \Rightarrow \nabla p$}
\LeftLabel{\small(w)'s}
\UnaryInfC{$\neg(p\wedge q), \nabla p \Rightarrow \nabla p, p\wedge q, \neg\nabla q$}
\LeftLabel{\small($\Rightarrow\neg\nabla$)}
\UnaryInfC{$\neg(p\wedge q) \Rightarrow \neg\nabla p, \nabla p, p\wedge q, \neg\nabla q$}

\AxiomC{$\nabla q \Rightarrow \nabla q$}
\LeftLabel{\small(w)'s}
\UnaryInfC{$\neg(p\wedge q), \nabla q \Rightarrow \nabla q, p\wedge q, \neg\nabla p$}
\LeftLabel{\small($\Rightarrow\neg\nabla$)}
\UnaryInfC{$\neg(p\wedge q) \Rightarrow  \neg\nabla q, \nabla q,  p\wedge q, \neg\nabla p$}

\LeftLabel{\small($\Rightarrow\nabla\wedge$)}
\BinaryInfC{$\neg(p\wedge q) \Rightarrow \nabla (p\wedge q), \neg\nabla p, p\wedge q, \neg\nabla q$}
\LeftLabel{\small($\nabla\neg\nabla\Rightarrow$)}
\UnaryInfC{$\nabla\neg\nabla(p\wedge q) \Rightarrow  \neg\nabla p, \neg\nabla q$}
\LeftLabel{\small($\Rightarrow\vee$)}
\UnaryInfC{$\nabla\neg\nabla(p\wedge q) \Rightarrow  \neg\nabla p \vee \neg\nabla q$}
\end{prooftree}
\end{eje}

\noindent It is clear that for any sequent $\Gamma\Rightarrow\Delta$, $\gsub{\Gamma\Rightarrow\Delta}$ is a finite set.\\[2mm]
\SetKwFunction{LowSeqGSix}{LowSeqGSix}
Let $\cal S$ and ${\cal S}_j$ be two sets of sequents. Let us denote by \LowSeqGSix(${\cal S}$, ${\cal S}_j)$ the set of all sequents $\Pi\Rightarrow\Lambda$ in  ${\cal S} - {\cal S}_{j}$ such that $\Pi\Rightarrow\Lambda$  is the lower sequent  of an inference in \GSix\  whose upper sequent(s) is (are) one (two) in ${\cal S}_{j}$. Then, consider the following algorithm written in pseudocode.

\

\IncMargin{1em}
\begin{algorithm}[H]
\SetKwInOut{Input}{input}\SetKwInOut{Output}{output}
\Input{ a sequent  $\Gamma \Rightarrow \Delta$ }
\Output{provable / not provable}
\BlankLine
${\cal S}  \longleftarrow \{ \Pi\Rightarrow\Lambda \, | \, \Pi\cup\Lambda \in \gsub{\Gamma \Rightarrow \Delta} \}$\;
${\cal S}_{0} \longleftarrow \{ \Pi\Rightarrow\Lambda \in {\cal S} \, | \,  \Pi\Rightarrow\Lambda \mbox{ is an initial sequent of \GSix}\}$\;
$j  \longleftarrow 0$\;
${\cal S}_{-1}  \longleftarrow \emptyset$\;
\While{$\Gamma \Rightarrow \Delta \notin {\cal S}_{j}$ and ${\cal S}_{j}\not= {\cal S}_{j-1}$}{
${\cal S}_{j+1}\longleftarrow \LowSeqGSix({\cal S}, {\cal S}_j) \cup {\cal S}_{j}$\;
$S\longleftarrow S - {\cal S}_{j}$\;
$j\longleftarrow j + 1$\;}
\lIf{$\Gamma \Rightarrow \Delta \in {\cal S}_{j}$}{\Return provable}
\lElse {\Return not provable}

\caption{Decision procedure for \GSix}\label{algo_disjdecomp}
\end{algorithm}\DecMargin{1em}

\

\noindent Since the set $\{ \Pi\Rightarrow\Lambda \, | \, \Pi\cup\Lambda \in \gsub{\Gamma \Rightarrow \Delta} \}$ is finite, this procedure always ends. 


\

\section{Conclusions}
We proposed to find a cut--free sequent system for the logic \Six\ using a method due to Avron, Ben-Naim and Konikowska \cite{Avron02} and taking advantage of the fact that \Six\ is determined by a single 6-element  logic matrix \cite{SMUR}. As a consequence of this process, we obtained a system with 84 logic rules (\SFSix ), a system with 230 logic rules (\TSSix ) and, after streamlining the latter, the system \GSix\ with 26 logic rules; all of them enjoying the cut--elimination property. The fact that we were dealing with six truth values increased significantly the complexity of the process, in particular, the streamlining process. In this respect, we propose (and prove their soundness) specific tools, such as Proposition \ref{propred} and Principle 4,  that can be used systematically to streamline the system. However, the success of this task deeply depends on the way one prove that language is sufficiently expressive. That is, it deeply depends on the choice we make for the $\alpha$'s and $\beta$'s formulas of Definition \ref{SufficientlyExp}. In particular, it is desirable to find  formulas  such that, for two given truth values, the respective formulas for $i$ and $j$ coincide except for one which in one case is an $\alpha$ and in the other is a $\beta$ (see the respective formulas for the values $\dt$ and $\1$ in Proposition \ref{AsBs}). 
After using our principles, one can use the semantics (along with Principle 2) to further streamline the system. Finally, we use the cut-free sequent system \GSix\ to provide a decision procedure for it and from which we obtain bottom-up proof search. We leave for a future work the study of the complexity of our method as well as its comparison with both, the decision procedure displayed in \cite{SMUR} and the one of truth-tables that is available from the matrix ${\cal M}_6$. 

We propose to continue our study by extending this method to logics that preserve degrees of truth w.r.t. some ordered structure but which do not have a matrix semantics. Also, we shall continue looking for tools that allow further systematization of the streamlining process.\\[6mm]
{\bf Acknowledgements.} We thank the anonymous referees for their careful reading, helpful suggestions and constructive comments
on this paper.

\end{document}